\documentclass[12pt]{amsart}
\usepackage{amssymb}
\begin{document}
\theoremstyle{plain}
\newtheorem{thm}{Theorem}[section]
\newtheorem{theorem}[thm]{Theorem}
\newtheorem{lemma}[thm]{Lemma}
\newtheorem{corollary}[thm]{Corollary}
\newtheorem{proposition}[thm]{Proposition}
\newtheorem{addendum}[thm]{Addendum}
\newtheorem{variant}[thm]{Variant}
\theoremstyle{definition}
\newtheorem{construction}[thm]{Construction}
\newtheorem{notations}[thm]{Notations}
\newtheorem{question}[thm]{Question}
\newtheorem{problem}[thm]{Problem}
\newtheorem{remark}[thm]{Remark}
\newtheorem{remarks}[thm]{Remarks}
\newtheorem{definition}[thm]{Definition}
\newtheorem{claim}[thm]{Claim}
\newtheorem{assumption}[thm]{Assumption}
\newtheorem{assumptions}[thm]{Assumptions}
\newtheorem{properties}[thm]{Properties}
\newtheorem{example}[thm]{Example}
\numberwithin{equation}{thm}
\catcode`\@=11
\def\opn#1#2{\def#1{\mathop{\kern0pt\fam0#2}\nolimits}}
\def\bold#1{{\bf #1}}%
\def\underrightarrow{\mathpalette\underrightarrow@}
\def\underrightarrow@#1#2{\vtop{\ialign{$##$\cr
 \hfil#1#2\hfil\cr\noalign{\nointerlineskip}%
 #1{-}\mkern-6mu\cleaders\hbox{$#1\mkern-2mu{-}\mkern-2mu$}\hfill
 \mkern-6mu{\to}\cr}}}
\let\underarrow\underrightarrow
\def\underleftarrow{\mathpalette\underleftarrow@}
\def\underleftarrow@#1#2{\vtop{\ialign{$##$\cr
 \hfil#1#2\hfil\cr\noalign{\nointerlineskip}#1{\leftarrow}\mkern-6mu
 \cleaders\hbox{$#1\mkern-2mu{-}\mkern-2mu$}\hfill
 \mkern-6mu{-}\cr}}}
\let\amp@rs@nd@\relax
\newdimen\ex@
\ex@.2326ex
\newdimen\bigaw@
\newdimen\minaw@
\minaw@16.08739\ex@
\newdimen\minCDaw@
\minCDaw@2.5pc
\newif\ifCD@
\def\minCDarrowwidth#1{\minCDaw@#1}
\newenvironment{CD}{\@CD}{\@endCD}
\def\@CD{\def\A##1A##2A{\llap{$\vcenter{\hbox
 {$\scriptstyle##1$}}$}\Big\uparrow\rlap{$\vcenter{\hbox{%
$\scriptstyle##2$}}$}&&}%
\def\V##1V##2V{\llap{$\vcenter{\hbox
 {$\scriptstyle##1$}}$}\Big\downarrow\rlap{$\vcenter{\hbox{%
$\scriptstyle##2$}}$}&&}%
\def\={&\hskip.5em\mathrel
 {\vbox{\hrule width\minCDaw@\vskip3\ex@\hrule width
 \minCDaw@}}\hskip.5em&}%
\def\verteq{\Big\Vert&&}%
\def\noarr{&&}%
\def\vspace##1{\noalign{\vskip##1\relax}}\relax\let\amp@rs@nd@&\iffalse}\fi
 \CD@true\vcenter\bgroup\relax\let\\=\cr\iffalse}\fi\tabskip\z@skip\baselineskip20\ex@
 \lineskip3\ex@\lineskiplimit3\ex@\halign\bgroup
 &\hfill$\m@th##$\hfill\cr}
\def\@endCD{\cr\egroup\egroup}
\def\>#1>#2>{\amp@rs@nd@\setbox\z@\hbox{$\scriptstyle
 \;{#1}\;\;$}\setbox\@ne\hbox{$\scriptstyle\;{#2}\;\;$}\setbox\tw@
 \hbox{$#2$}\ifCD@
 \global\bigaw@\minCDaw@\else\global\bigaw@\minaw@\fi
 \ifdim\wd\z@>\bigaw@\global\bigaw@\wd\z@\fi
 \ifdim\wd\@ne>\bigaw@\global\bigaw@\wd\@ne\fi
 \ifCD@\hskip.5em\fi
 \ifdim\wd\tw@>\z@
 \mathrel{\mathop{\hbox to\bigaw@{\rightarrowfill}}\limits^{#1}_{#2}}\else
 \mathrel{\mathop{\hbox to\bigaw@{\rightarrowfill}}\limits^{#1}}\fi
 \ifCD@\hskip.5em\fi\amp@rs@nd@}
\def\<#1<#2<{\amp@rs@nd@\setbox\z@\hbox{$\scriptstyle
 \;\;{#1}\;$}\setbox\@ne\hbox{$\scriptstyle\;\;{#2}\;$}\setbox\tw@
 \hbox{$#2$}\ifCD@
 \global\bigaw@\minCDaw@\else\global\bigaw@\minaw@\fi
 \ifdim\wd\z@>\bigaw@\global\bigaw@\wd\z@\fi
 \ifdim\wd\@ne>\bigaw@\global\bigaw@\wd\@ne\fi
 \ifCD@\hskip.5em\fi
 \ifdim\wd\tw@>\z@
 \mathrel{\mathop{\hbox to\bigaw@{\leftarrowfill}}\limits^{#1}_{#2}}\else
 \mathrel{\mathop{\hbox to\bigaw@{\leftarrowfill}}\limits^{#1}}\fi
 \ifCD@\hskip.5em\fi\amp@rs@nd@}
\newenvironment{CDS}{\@CDS}{\@endCDS}
\def\@CDS{\def\A##1A##2A{\llap{$\vcenter{\hbox
 {$\scriptstyle##1$}}$}\Big\uparrow\rlap{$\vcenter{\hbox{%
$\scriptstyle##2$}}$}&}%
\def\V##1V##2V{\llap{$\vcenter{\hbox
 {$\scriptstyle##1$}}$}\Big\downarrow\rlap{$\vcenter{\hbox{%
$\scriptstyle##2$}}$}&}%
\def\={&\hskip.5em\mathrel
 {\vbox{\hrule width\minCDaw@\vskip3\ex@\hrule width
 \minCDaw@}}\hskip.5em&}
\def\verteq{\Big\Vert&}
\def\novarr{&}
\def\noharr{&&}
\def\SE##1E##2E{\slantedarrow(0,18)(4,-3){##1}{##2}&}
\def\SW##1W##2W{\slantedarrow(24,18)(-4,-3){##1}{##2}&}
\def\NE##1E##2E{\slantedarrow(0,0)(4,3){##1}{##2}&}
\def\NW##1W##2W{\slantedarrow(24,0)(-4,3){##1}{##2}&}
\def\slantedarrow(##1)(##2)##3##4{%
\thinlines\unitlength1pt\lower 6.5pt\hbox{\begin{picture}(24,18)%
\put(##1){\vector(##2){24}}%
\put(0,8){$\scriptstyle##3$}%
\put(20,8){$\scriptstyle##4$}%
\end{picture}}}
\def\vspace##1{\noalign{\vskip##1\relax}}\relax\let\amp@rs@nd@&\iffalse}\fi
 \CD@true\vcenter\bgroup\relax\let\\=\cr\iffalse}\fi\tabskip\z@skip\baselineskip20\ex@
 \lineskip3\ex@\lineskiplimit3\ex@\halign\bgroup
 &\hfill$\m@th##$\hfill\cr}
\def\@endCDS{\cr\egroup\egroup}
\newdimen\TriCDarrw@
\newif\ifTriV@
\newenvironment{TriCDV}{\@TriCDV}{\@endTriCD}
\newenvironment{TriCDA}{\@TriCDA}{\@endTriCD}
\def\@TriCDV{\TriV@true\def\TriCDpos@{6}\@TriCD}
\def\@TriCDA{\TriV@false\def\TriCDpos@{10}\@TriCD}
\def\@TriCD#1#2#3#4#5#6{%
\setbox0\hbox{$\ifTriV@#6\else#1\fi$}
\TriCDarrw@=\wd0 \advance\TriCDarrw@ 24pt
\advance\TriCDarrw@ -1em
\def\SE##1E##2E{\slantedarrow(0,18)(2,-3){##1}{##2}&}
\def\SW##1W##2W{\slantedarrow(12,18)(-2,-3){##1}{##2}&}
\def\NE##1E##2E{\slantedarrow(0,0)(2,3){##1}{##2}&}
\def\NW##1W##2W{\slantedarrow(12,0)(-2,3){##1}{##2}&}
\def\slantedarrow(##1)(##2)##3##4{\thinlines\unitlength1pt
\lower 6.5pt\hbox{\begin{picture}(12,18)%
\put(##1){\vector(##2){12}}%
\put(-4,\TriCDpos@){$\scriptstyle##3$}%
\put(12,\TriCDpos@){$\scriptstyle##4$}%
\end{picture}}}
\def\={\mathrel {\vbox{\hrule
   width\TriCDarrw@\vskip3\ex@\hrule width
   \TriCDarrw@}}}
\def\>##1>>{\setbox\z@\hbox{$\scriptstyle
 \;{##1}\;\;$}\global\bigaw@\TriCDarrw@
 \ifdim\wd\z@>\bigaw@\global\bigaw@\wd\z@\fi
 \hskip.5em
 \mathrel{\mathop{\hbox to \TriCDarrw@
{\rightarrowfill}}\limits^{##1}}
 \hskip.5em}
\def\<##1<<{\setbox\z@\hbox{$\scriptstyle
 \;{##1}\;\;$}\global\bigaw@\TriCDarrw@
 \ifdim\wd\z@>\bigaw@\global\bigaw@\wd\z@\fi
 \mathrel{\mathop{\hbox to\bigaw@{\leftarrowfill}}\limits^{##1}}
 }
 \CD@true\vcenter\bgroup\relax\let\\=\cr\iffalse}\fi
 \tabskip\z@skip\baselineskip20\ex@
 \lineskip3\ex@\lineskiplimit3\ex@
 \ifTriV@
 \halign\bgroup
 &\hfill$\m@th##$\hfill\cr
#1&\multispan3\hfill$#2$\hfill&#3\\
&#4&#5\\
&&#6\cr\egroup%
\else
 \halign\bgroup
 &\hfill$\m@th##$\hfill\cr
&&#1\\%
&#2&#3\\
#4&\multispan3\hfill$#5$\hfill&#6\cr\egroup
\fi}
\def\@endTriCD{\egroup}
\newcommand{\sA}{{\mathcal A}}
\newcommand{\sB}{{\mathcal B}}
\newcommand{\sC}{{\mathcal C}}
\newcommand{\sD}{{\mathcal D}}
\newcommand{\sE}{{\mathcal E}}
\newcommand{\sF}{{\mathcal F}}
\newcommand{\sG}{{\mathcal G}}
\newcommand{\sH}{{\mathcal H}}
\newcommand{\sI}{{\mathcal I}}
\newcommand{\sJ}{{\mathcal J}}
\newcommand{\sK}{{\mathcal K}}
\newcommand{\sL}{{\mathcal L}}
\newcommand{\sM}{{\mathcal M}}
\newcommand{\sN}{{\mathcal N}}
\newcommand{\sO}{{\mathcal O}}
\newcommand{\sP}{{\mathcal P}}
\newcommand{\sQ}{{\mathcal Q}}
\newcommand{\sR}{{\mathcal R}}
\newcommand{\sS}{{\mathcal S}}
\newcommand{\sT}{{\mathcal T}}
\newcommand{\sU}{{\mathcal U}}
\newcommand{\sV}{{\mathcal V}}
\newcommand{\sW}{{\mathcal W}}
\newcommand{\sX}{{\mathcal X}}
\newcommand{\sY}{{\mathcal Y}}
\newcommand{\sZ}{{\mathcal Z}}
\newcommand{\A}{{\mathbb A}}
\newcommand{\B}{{\mathbb B}}
\newcommand{\C}{{\mathbb C}}
\newcommand{\D}{{\mathbb D}}
\newcommand{\E}{{\mathbb E}}
\newcommand{\F}{{\mathbb F}}
\newcommand{\G}{{\mathbb G}}
\newcommand{\HH}{{\mathbb H}}
\newcommand{\I}{{\mathbb I}}
\newcommand{\J}{{\mathbb J}}
\renewcommand{\L}{{\mathbb L}}
\newcommand{\M}{{\mathbb M}}
\newcommand{\N}{{\mathbb N}}
\renewcommand{\P}{{\mathbb P}}
\newcommand{\Q}{{\mathbb Q}}
\newcommand{\R}{{\mathbb R}}
\newcommand{\T}{{\mathbb T}}
\newcommand{\U}{{\mathbb U}}
\newcommand{\V}{{\mathbb V}}
\newcommand{\W}{{\mathbb W}}
\newcommand{\X}{{\mathbb X}}
\newcommand{\Y}{{\mathbb Y}}
\newcommand{\Z}{{\mathbb Z}}
\newcommand{\id}{{\rm id}}
\newcommand{\rank}{{\rm rank}}
\newcommand{\END}{{\mathbb E}{\rm nd}}
\newcommand{\End}{{\rm End}}
\newcommand{\Hg}{{\rm Hg}}
\newcommand{\tr}{{\rm tr}}
\newcommand{\Sl}{{\rm Sl}}
\newcommand{\Gl}{{\rm Gl}}
\newcommand{\Cor}{{\rm Cor}}
\title[A characterization of Shimura curves]{A characterization of
certain Shimura curves in the moduli stack of abelian varieties}
\author[Eckart Viehweg]{Eckart Viehweg}
\address{Universit\"at Essen, FB6 Mathematik, 45117 Essen, Germany}
\email{ viehweg@uni-essen.de}
\thanks{This work has been supported by the ``DFG-Schwerpunktprogramm
Globale Methoden in der Komplexen Geometrie''. The second named author
is supported by a grant from the Research
Grants Council of the Hong Kong
Special Administrative Region, China
(Project No. CUHK 4239/01P)}
\author[Kang Zuo]{Kang Zuo}
\address{The Chinese University of Hong Kong, Department of Mathematics,
Shatin, Hong Kong}
\email{kzuo@math.cuhk.edu.hk}
\maketitle
Throughout this article, $Y$ will denote a non-singular complex projective
curve, $U$ an open dense subset, and $X_0\to U$ a smooth family of abelian varieties.
We choose a projective non-singular compactification $X$ of $X_0$ such that the
family extends to a morphism $f:X\to Y$, which we call again a family
of abelian varieties although some of the fibres are singular. We write
$S = Y\setminus U$, and $\Delta=f^{-1}(S)$. Consider the
weight $1$ variation of Hodge structures given by $f:X_0\to U$,
i.e. $R^1f_*\Z_{X_0}$. We will always assume that the monodromy
of $R^1f_*\Z_{X_0}$ around all points in $S$ is unipotent. The
Deligne extension of $(R^1f_*\Z_{X_0})\otimes \sO_U$ to $Y$
carries a Hodge filtration. Taking the graded sheaf one obtains
the Higgs bundle
$$
(E, \theta)=(E^{1,0}\oplus E^{0,1}, \theta_{1,0})
$$
with $E^{1,0}=f_*\Omega^1_{X/Y}(\log \Delta)$ and
$E^{0,1}=R^1f_*\sO_X$.
The Higgs field $\theta_{1,0}$ is given by the edge morphisms
$$
f_*\Omega^1_{X/Y}(\log \Delta)\>>> R^1f_*\sO_X \otimes
\Omega^{1}_{Y}(\log S)
$$
of the tautological sequence
$$
0\to {f}^*\Omega^1_Y(\log S)\to \Omega^1_{X}(\log \Delta) \to
\Omega_{X/Y}^1(\log \Delta))\to 0.
$$
By \cite{Kol} $E$ can be decomposed as a direct sum $F\oplus N$ of Higgs bundles
with $E^{1,0}\cap F$ ample and with $N$ flat,
hence for $F^{i,j}=E^{i,j}\cap F$ and $N^{i,j}=E^{i,j}\cap N$
the Higgs bundle $E$ decomposes in
\begin{equation}\label{split1}
(F=F^{1,0}\oplus F^{0,1}, \theta_{1,0}|_{F^{1,0}}) \mbox{ \ \ \
and \ \ \ } (N^{1,0}\oplus N^{0,1},0).
\end{equation}

For $g_0={\rm rank}(F^{1,0})$ the Arakelov inequalities
(\cite{De4}, generalized in \cite{Pet}, \cite{J-Z}) say that
\begin{equation}\label{arakelovineq}
2\cdot \deg(F^{1,0}) \leq g_0\cdot \deg(\Omega^1_Y(\log S)).
\end{equation}
In this note we will try to understand $f:X\to Y$, for which
(\ref{arakelovineq}) is an equality, or as we will say, of
families reaching the Arakelov bound. By \ref{semistablei}, this
property is equivalent to the maximality of the Higgs field for
$F$, saying that $\theta_{1,0}:F^{1,0} \to F^{0,1}\otimes
\Omega^1_Y(\log S)$ is an isomorphism.

As it will turn out, the base of a family of abelian varieties
reaching the Arakelov bound is a Shimura curve, and the
maximality of the Higgs field is reflected in the existence of special
Hodge cycles on the general fibre. Before formulating a general
result, let us consider two examples.

For families of elliptic curves, the maximality of the Higgs field
just says that the family is modular (see Section \ref{sectshimura}).

\begin{proposition}\label{modular}
Let $f:E \to Y$ be a semi-stable family of elliptic curves, smooth
over $U\subset Y$. If $E\to Y$ is non isotrivial and reaching the
Arakelov bound, $E\to Y$ is modular, i.e. $U$ is the quotient of
the upper half plane $\sH$ by a subgroup of $\text{Sl}_2(\Z)$ of
finite index, and the morphism $U \to \C=\sH/ \text{Sl}_2(\Z)$ is
given by the $j$-invariant of the fibres.
\end{proposition}

For $S\neq \emptyset$ the only families of abelian varieties
reaching the Arakelov bound are build up from modular families of
elliptic curves.

\begin{theorem}\label{geomsplit}
Let $f:X\to Y$ be a family of abelian varieties smooth over $U$,
and such that the local monodromies around $s\in S$ are unipotent.
If $S\not=\emptyset,$ and if $f:X\to Y$
reaches the Arakelov bound, then there exists an \'etale covering
$\pi:Y'\to Y$ such that $f':X'=X\times_YY'\to Y'$ is isogenous over $Y'$
to a product
$$ B\times E\times_{Y'}\cdots\times_{Y'} E ,$$
where $B$ is abelian variety defined over $\C$ of dimension
$g-g_0,$ and where $h: E\to Y'$ is a family of semi-stable elliptic
curves reaching the Arakelov bound.
\end{theorem}

Results parallel to \ref{geomsplit} have been obtained in
\cite{STZ} for families of $K3$-surfaces, and the methods and
results of \cite{STZ} have been a motivation to study the case of abelian
varieties.

As we will see in Section \ref{sectqsplitting} Theorem
\ref{geomsplit} follows from the existence of too many
endomorphisms of the general fibre of $f:X\to Y$, which in turn
implies the existence of too many cycles on the
general fibre of $X\times_Y X$. We give an elementary
proof of Theorem \ref{geomsplit} in Section \ref{sectqsplitting},
although it is nothing but
a first example for the relation between the maximality
of Higgs fields, and the moduli of abelian varieties with a
given special Mumford-Tate group $\Hg$, constructed in \cite{Mum1} and \cite{Mum2}
(see Section \ref{sectshimura}).

\begin{proposition}\label{MTfamilies}
Let $f: X\to Y$ be a family of $g$-dimensional abelian varieties
reaching the Arakelov bound. Assume that $g=g_0$, or more generally that the
largest unitary local subsystem $\U_1$ of $R^1f_*\C_{X_0}$ is defined over $\Q$.
Then there exists a finite cover $Y'\to Y$, \'etale over $U$, and a $\Q$-algebraic
subgroup $\Hg \subset {\rm Sp}(2g,\R)$, such that pullback family
$f': Y'=X\times_YY' \to Y'$ is a semi-stable compactification of the universal family of
polarized abelian varieties with special Mumford-Tate group contained in $\Hg$,
and with a suitable level structure.
\end{proposition}

As a preparation for the proof of Proposition \ref{MTfamilies}
we will show in Section \ref{sectsplitting}, using
Simpson's correspondence between Higgs bundles and local systems,
that the maximality of the Higgs field enforces a presentation of the
local systems $R^1f_*\C_{X_0}$ and $\END(R^1f_*\C_{X_0})$ using
direct sums and tensor products of one weight one
complex variation of Hodge structures $\L$ of rank two and several
unitary local systems.

Proposition \ref{MTfamilies} relates families 
reaching the Arakelov bound to totally geodesic subvarieties of the moduli 
space of abelian varieties, as considered by Moonen in \cite{Mo}, or to the 
totally geodesic holomorphic embeddings, studied by Abdulali in \cite{Abd}
(see Remark \ref{remark}, b). As in \cite{Sai}
one could use the classification of Shimura varieties due to Satake
\cite{Sat} to obtain a complete list of those families, and to characterize
them in terms of properties of their variation of Hodge structures.

We choose a different approach, less relying on the theory of Shimura varieties,
and more adapted to handle the remaining families of abelian varieties
(see Remark \ref{remark}, c),
as well as some other families of varieties of Kodaira dimension zero (see \cite{VZ}).
We first show that the decompositions of $R^1f_*\C_{X_0}$ and $\END(R^1f_*\C_{X_0})$
mentioned above are defined over $\bar\Q\cap\R$.
In case $S\neq \emptyset$ it is then easy to see, that the
unitary parts of the decompositions trivialize, after replacing
$Y$ by a finite \'etale cover $Y'$ (see \ref{split9}).

For $S=\emptyset$, let us assume first that the assumptions made in \ref{MTfamilies}
hold true. By \cite{Fal} (see Proposition \ref{rigidity}) they imply that the
family is rigid, i.e. that the morphism from $Y$ to the moduli
stack of polarized abelian varieties has no non-trivial
deformation, except those obtained by deforming a constant abelian
subvariety.

Mumford gave in \cite{Mum2} countably many moduli functors
of abelian fourfolds, where $\Hg$ is obtained via the corestriction
of an quaternion algebra, defined over a totally real cubic number
field $F$. Generalizing his construction one considers
quaternion division algebras $A$ defined over any totally real
number field $F$, which are ramified at all infinite places
except one. Choose an embedding
$$
D=\Cor_{F/\Q}A \subset M(2^m,\Q),
$$
with $m$ minimal. As we will see in Section \ref{sectquaternion}
writing $d=[F:\Q]$ one finds $m=d$ or $m=d+1$. By \ref{constrq1}
and \ref{constrq2} we get the following types of moduli functors
of abelian varieties with special Mumford-Tate group
$$
\Hg=\{ x\in D^*; \ x{\bar x}=1\}
$$
and with a suitable level structure, which are represented by a
smooth family $Z_A \to Y_A $ over a compact Shimura curve $Y_A$.
Since we did not fix the level structure, $Y_A$ is not uniquely
determined by $A$. So it rather stands for a whole class of
possible base curves, two of which have a common finite \'etale
covering.
\begin{example}\label{mumfordshimura} \ Let $Z_\eta$ denote the generic
fibre of $Z_A \to Y_A$. Then one of the following holds true.
\begin{enumerate}
\item[i.] $1< m=d$ odd. In this case
$\dim(Z_\eta)=2^{d-1}$ and ${\rm End}(Z_\eta)\otimes_\Z \Q = \Q$.
\item[ii.] $m=d+1$. Then $\dim(Z_\eta)=2^{d}$ and
\begin{enumerate}
\item[a.] for $d$ odd, ${\rm End}(Z_\eta)
\otimes_\Z\Q$ a totally indefinite quaternion algebra over $\Q$.
\item[b.] for $d$ even, ${\rm End}(Z_\eta)
\otimes_\Z\Q$ a totally definite quaternion algebra over $\Q$.
\end{enumerate}
\end{enumerate}
\end{example}
Let us call the family $Z_A\to Y_A$ a family of Mumford type.

For $d=1$ or $2$ the examples in \ref{mumfordshimura} include the only two Shimura curves
of PEL-type, parameterizing
\begin{enumerate}
\item[$\bullet$] Moduli schemes of false elliptic curves, i.e. polarized abelian
surfaces $B$ with ${\rm End}(B)\otimes_\Z\Q$ a totally indefinite
quaternion algebra over $\Q$ (see also \cite{Sha}).
\item[$\bullet$] Moduli schemes of abelian fourfolds $B$ with ${\rm End}(B)
\otimes_\Z\Q$ a totally definite quaternion algebra over $\Q$.
\end{enumerate}

We will see in Section \ref{sectcorestriction}
that for $g=g_0$, up to powers and isogenies,
the families of Mumford type are the only smooth families of abelian
varieties over curves reaching the Arakelov bound.

\begin{theorem}\label{shimura}
Let $f:X\to Y$ be a smooth family of abelian varieties. If the
largest unitary local subsystem $\U_1$ of $R^1f_*\C_X$ is
defined over $\Q$ and if $f:X\to Y$ reaches the Arakelov bound,
then there exist
\begin{enumerate}
\item[a.] a quaternion division algebra $A$, defined over a totally
real number field $F$, and ramified at all infinite places except one,
\item[b.] an \'etale covering $\pi:Y'\to Y$,
\item[c.] a family of Mumford type $h:Z=Z_A \to Y'=Y_A$, as in Example
\ref{mumfordshimura}, and an abelian variety $B$ such that
$f':X'=X\times_YY'\to Y'$ is isogenous to
$$
B\times Z\times_{Y'} \cdots \times_{Y'}Z \>>> Y'.
$$
\end{enumerate}
\end{theorem}

Things are getting more complicated if one drops the condition on
the unitary local subsystem $\U_1$ of $R^1f_*\C_X$. For one quaternion algebra $A$,
there exist several non isogenous families. Hence it will no longer be true, that up to a
constant factor $f:X\to Y$ is isogenous to the product of one particular family.

\begin{example}[see \ref{constrl}]\label{general}
Let $A$ be a quaternion algebra defined over a totally real
number field $F$, ramified at all infinite places but one, and let
$L$ be a subfield of $F$. Let $\beta_1,\ldots,\beta_\delta:L\to
\bar\Q$ denote the different embeddings of $L$. For $\mu=[F:L]+1$
(or may be $\mu=[F:\Q]$ in case that $L=\Q$ and $\mu$ odd) there
exists an embedding
$$
\Cor_{F/L}A \subset M(2^\mu,L).
$$
As well known (see Section \ref{sectquaternion}) for some
Shimura curve $Y'$ such an embedding gives rise to a
representation of $\pi_1(Y',*)$ in $M(2^\mu,L)$, hence to a local
$L$ system $\V_L$. Moreover there exists an irreducible $\Q$ local system
$\X_\Q=\X_{A,L;\Q}$ for which $\X_\Q\otimes \bar\Q$ is a direct
sum of the local systems $\V_L\otimes_{L,\beta_\nu}\bar\Q$.
\end{example}
There exist non-isotrivial families $h:Z\to Y'$ with a geometrically simple
generic fibre, such that $R^1h_*\Q$ is a direct sum of $\iota$
copies of $\X_\Q$. Such examples, for $g=4$ and $8$ have been
considered in \cite{Fal}. Here $F$ is a quadratic extension of
$\Q$, $L=F$ and $\iota=1$ or $2$. For $g=8$, i.e. $\iota=2$, this
gives the lowest dimensional example of a non rigid family of
abelian varieties without a trivial sub family \cite{Fal}. A
complete classification of such families is given in \cite{Sai}.

\begin{theorem}\label{shimura2}
Let $f:X\to Y$ be a smooth family of abelian varieties. If $f:X\to
Y$ reaches the Arakelov bound, then there exist an \'etale
covering $\pi:Y'\to Y$, a quaternion algebra $A$, defined over a
totally real number field $F$ and ramified at all of the infinite
places except one, an abelian variety $B$, and $\ell$ families
$h_i:Z_i \to Y'$ of abelian varieties with geometrically simple generic fibre,
such that
\begin{enumerate}
\item[i.] $f':X'=X\times_YY'\to Y'$ is isogenous to
$$
B\times Z_1 \times_{Y'} \cdots \times_{Y'}Z_\ell \>>> Y'.
$$
\item[ii.] For each $i\in\{1,\ldots,\ell\}$ there exists a subfield
$L_i$ of $F$ such that the local system $R^1h_{i*}\Q_{Z_i}$ is a direct
sum of copies of the irreducible $\Q$ local system $\X_{A,L_i;\Q}$
defined in Example \ref{general}.
\item[iii.] For each $i\in\{1, \ldots, \ell\}$ the following conditions
are equivalent:
\begin{enumerate}
\item[a.] $L_i=\Q$.
\item[b.] $h_i:Z_i \to Y'$ is a family of Mumford type, as defined
in Example \ref{mumfordshimura}.
\item[c.] $\End(\X_{A,L_i;\Q})^{0,0}=\End(\X_{A,L_i;\Q})$.
\end{enumerate}
Moreover, if one of those conditions holds true, $R^1h_{i*}\Q_{Z_i}$
is irreducible, hence $R^1h_{i*}\Q_{Z_i}=\X_{A,L_i;\Q}$.
\end{enumerate}
\end{theorem}
Here, contrary to \ref{shimura}, we do not claim that a component
$h_i:Z_i \to Y'$ is uniquely determined up to isogeny by
$\X_{A,L_i;\Q}$ and by the rank of $R^1h_{i*}\Q_{Z_i}$.\\

We do not know for which $g$ there are families of
Jacobians among the families of abelian varieties considered in
\ref{geomsplit}, \ref{shimura} or \ref{shimura2}, i.e. whether
one can find a family $\varphi:Z\to Y$ of curves of genus $g$ such that
$f:J(Z/Y)\to Y$ reaches the Arakelov bound.

For $Y=\P^1$ the Arakelov inequality (\ref{arakelovineq}) implies
that $\#S\geq 4$. Our hope, that a family of abelian varieties
with $\# S=4$ can not be a family of Jacobians, broke down when we
found an example of a family of genus $2$ curves over the
modular curve $X(3)$ in \cite{Kan1}, whose Jacobian is isogenous to the
product of a fixed elliptic curve $B$ with the modular curve
$E(3)\to X(3)$ (see Section \ref{sectjacobians}).\\

As mentioned already, this article owes a lot to the recent work of
the second named author with Xiao-Tao Sun and Sheng-Li Tan.
We thank Ernst Kani for explaining his beautiful construction in \cite{Kan1},
and for sharing his view about higher genus analogs of families of curves with
splitting Jacobians. It is also a pleasure to thank Ben Moonen, H\'el\`ene Esnault
and Frans Oort for their interest and help, Ngaiming Mok, for explaining
us differential geometric properties of base spaces of families and for pointing
out Mumford's construction in \cite{Mum2}, Bruno Kahn and Claus Scheiderer
for their help to understand quaternion algebras and their corestriction.

This note grew out of discussions started when the first named author visited
the Institute of Mathematical Science and the Department of Mathematics at the
Chinese University of Hong Kong. His contributions to the final version
(in particular to the proof of Theorems \ref{shimura} and \ref{shimura2})
were written during his visit to the I.H.E.S., Bures sur Yvette.
He would like to thank the members of those three Institutes for their
hospitality.

\section{Splitting of $\C$-local systems} \label{sectsplitting}

We will frequently use C. Simpson's correspondence between
poly-stable Higgs bundles of degree zero and representations of
the fundamental group $\pi_1(U,*)$.

\begin{theorem}[C. Simpson \cite{Si2}]\label{simpson1}
There exists a natural equivalence between the category of direct
sums of stable filtered regular Higgs bundles of degree zero, and
of direct sums of stable filtered local systems of degree zero.
\end{theorem}
We will not recall the definition of a ``filtered regular'' Higgs bundle
(\cite{Si2}, p. 717), just remark that for a Higgs bundle corresponding
to a local system $\V$ with unipotent monodromy around the points in
$S$ the filtration is trivial, and automatically $\deg(\V)=0$.
By (\cite{Si2}, p. 720) the latter also holds true for
local systems $\V$ which are polarisable $\C$-variations of Hodge structures.

For example, \ref{simpson1} implies that the splitting of Higgs
bundles (\ref{split1}) corresponds to a decomposition over
$\mathbb C$
$$ (R^1f_*\Z_{X_0})\otimes \C=\V \oplus\U_1$$
where $\V$ corresponds to the Higgs bundle $(F=F^{1,0}\oplus
F^{0,1}, \theta)$ and $\U_1$ to $(N=N^{1,0}\oplus
N^{0,1},\theta_N=0)$. Let $\Theta(N,h)$  denote the curvature of
the Hodge metric $h$ on $E^{1,0}\oplus E^{0,1}$ restricted to $N,$
then by \cite{G}, chapter II we have
$$
\Theta(N,h|_N)=-\theta_N\wedge\bar\theta_N-\bar\theta_N\wedge\theta_N=0.
$$
This means that $h|_N$ is a flat metric. Hence,  $\U_1$ is a
unitary local system.

In general, if $\U$ is a local system, whose Higgs bundle is
a direct sum of stable Higgs bundles of degree zero and with a trivial
Higgs field, then $\U$ is unitary.

As a typical application of Simpson's correspondence one obtains the
polystability of the components of certain Higgs bundles. We just formulate
it in the weight one case.

Recall that $F^{1,0}$ is polystable, if there exists a
decomposition
$$ F^{1,0}\simeq \bigoplus_i \sA_i$$
with $\sA_i$ stable, and
$$\frac{\deg \sA_i}{\rank\sA_i}=\frac{\deg F^{1,0}}{\rank F^{1,0}}.$$

\begin{proposition}\label{semistablei}
Let $\V$ be a direct sum of stable filtered local systems of degree zero with
Higgs bundle $(F=F^{1,0}\oplus F^{0,1}, \tau)$. Assume that
$\tau|_{F^{0,1}}=0$, that
$$
\tau_{1,0}=\tau|_{F^{1,0}}:F^{1,0} \>>>
F^{0,1}\otimes\Omega^1_Y(\log S) \subset F\otimes\Omega^1_Y(\log
S),
$$
and that
\begin{equation}\label{arakeloveqC}
2 \cdot \deg(F^{1,0}) = g_0 \cdot \deg(\Omega^1_{Y}(\log S)).
\end{equation}
Then $\tau_{1,0}$ is an isomorphism, and the sheaf $F^{1,0}$ is
poly-stable.
\end{proposition}
\begin{proof}[Proof of \ref{semistablei}] Let $\sA\subset
F^{1,0}$ be a subsheaf, and let $\sB\otimes\Omega_Y^1(\log S)$ be
its image under $\theta_{1,0}$. Then $\sA\oplus\sB$ is a Higgs
subbundle of $F^{1,0}\oplus F^{0,1}$, and applying \ref{simpson1}
one finds $\deg(\sA)+\deg(\sB)\leq 0$. Hence
\begin{multline*}
\deg(\sA) = \deg(\sB) + \rank(\sB)\cdot\deg(\Omega_Y^1(\log S))\\
\leq \deg(\sB) + \rank(\sA)\cdot\deg(\Omega_Y^1(\log S)) \leq
-\deg(\sA) + \rank(\sA)\cdot\deg(\Omega_Y^1(\log S)).
\end{multline*}
The equality (\ref{arakeloveqC}) implies that
$$
\frac{\deg(\sA)}{\rank(\sA)} \leq \frac{1}{2}\deg(\Omega_Y^1(\log
S))= \frac{\deg(F^{1,0})}{g_0},
$$
and $F^{1,0}$ is semi-stable. If
$$
\frac{\deg(\sA)}{\rank(\sA)} = \frac{\deg(F^{1,0})}{g_0},
$$
${\rm rank}(\sA)={\rm rank}(\sB)$ and $\deg(\sB)=-\deg(\sA)$. The
Higgs bundle $(F^{1,0}\oplus F^{0,1},\theta)$ splits by
\ref{simpson1} as a direct sum of stable Higgs bundles of degree
zero. Hence $(\sA\oplus\sB,\theta|_{\sA\oplus\sB})$ is a direct
factor of $(F^{1,0}\oplus F^{0,1},\theta)$. In particular, $\sA$
is a direct factor of $F^{1,0}$. For $\sA=F^{1,0}$ one also
obtains that $\tau^{1,0}$ is injective and by (\ref{arakeloveqC})
it must be an isomorphism.
\end{proof}

The local system $R^1f_*\Q_{X_0}$ on $U=Y\setminus S$ is a $\Q$
variation of Hodge structures with unipotent local monodromies
around $s\in S$, obviously having a $\Z$-form.
By Deligne's semi-simplicity theorem \cite{Del}
it decomposes as a direct sum of
irreducible polarisable $\Q$-variation of Hodge structures $\mathbb V_{i\Q}$.

More generally, if $\V$ is a polarized $\C$-variation of Hodge structures, and
$$
\V =\bigoplus_i \mathbb V_{i},
$$
a decomposition with $\V_i$ an irreducible $\C$-local system, then by \cite{De4} each $\V_i$
again is a polarisable $\C$-variation of Hodge structures.

In both cases, taking the grading of the Hodge filtration, one obtains a
decomposition of the Higgs bundle
$$
(E,\theta)=(F^{1,0}\oplus F^{0,1}, \theta)\oplus (N^{1,0}\oplus
N^{0,1},0)
$$
as a direct sum of sub Higgs bundles, as stated in \ref{simpson1}.
Obviously, each of the $\V_{i\Q}$ again reaches the Arakelov bound. \\

Our next constructions will not require the local system to be
defined over $\Q$. So by abuse of notations, we will make the
following assumptions.

\begin{assumption}\label{assumption}
For a number field $L\subset \C$ consider a polarized $L$
variation of Hodge structures $\X_L$ of weight one over
$U=Y\setminus S$ with unipotent local monodromies around $s\in
S$. Assume that the local system $\X=\X_L\otimes_L\C$ has a
decomposition $\X=\V\oplus\U_1$, with $\U_1$ unitary,
corresponding to the decomposition
$$
(E,\theta)=(F,\theta_{1,0})\oplus (N,0)=(F^{1,0}\oplus F^{0,1},
\theta_{1,0}) \oplus (N^{1,0}\oplus N^{0,1},0)
$$
of Higgs fields. Assume that $\V$ (or $(F,\theta_{1,0})$) has a
maximal Higgs field, i.e. that
$$
\theta_{1,0}:F^{1,0} \to F^{0,1}\otimes\Omega^1_Y(\log S)
$$
is an isomorphism. Obviously, for $g_0={\rm rank}(F^{1,0})$ this
is equivalent to the equality (\ref{arakeloveqC}). Hence we will
also say, that $\X$ (or $(E,\theta)$) reaches the Arakelov bound.
\end{assumption}

\begin{proposition}\label{semistableii}
If $\deg\Omega^1_Y(\log S)$ is even there exists a tensor product
decomposition of variations of Hodge structures
$$ \V\simeq \L \otimes_\C \T,$$ with:
\begin{enumerate}
\item[a.] $\L$ is a rank-2 local system. For some invertible sheaf $\sL$, with
$\sL^2=\Omega^1_Y(\log S)$ the Higgs bundle corresponding to $\L$ is
$(\sL\oplus\sL^{-1},\tau)$, with $\tau|_{\sL^{-1}}=0$ and
$\tau|_{\sL}$ given by an isomorphism
$$
\tau^{1,0}: \sL\>>> \sL^{-1}\otimes\Omega^1_Y(\log S).
$$
$\sL$ has bidegree $1,0$, and $\sL^{-1}$ has bidegree $0,1$.
\item[b.] If $g_0$ is odd, $\sL^{g_0}=\det(F^{1,0})$ and $\sL$ is uniquely determined.
\item[c.] For $g_0$ even, there exists some invertible sheaf $\sN$ of order two in
${\rm Pic}^0(Y)$ with $\sL^{g_0} =\det(F^{1,0})\otimes \sN$.
\item[d.] $\T$ is a unitary local system and a variation of Hodge structures
of pure bidegree $0,0$. If $(\sT,0)$ denotes the
corresponding Higgs field, then $\sT=F^{1,0}\otimes\sL^{-1}=F^{0,1}\otimes\sL$.
\end{enumerate}
\end{proposition}
In section \ref{sectcorestriction} we will need a slightly
stronger statement.
\begin{addendum}\label{addendum}
If in \ref{semistableii}, there exists a presentation
$\V=\T_1\otimes_\C \V_1$ with $\T_1$ unitary and a variation of
Hodge structures of pure bidegree $0,0$, then there exists a
unitary local system $\T_2$ with $\T=\T_1\otimes_\C \T_2$.
\end{addendum}
In fact, write $(\sT_1,0)$ and $(F^{1,0}_1\oplus F^{0,1}_1,\theta_1)$ for
Higgs fields corresponding to $\T_1$ and $\V_1$, respectively. Then
$\deg(\sT_1)=0$ and
\begin{gather*}
2\cdot \deg(F_1^{1,0})\cdot {\rm rank}(\sT)=
2 \cdot \deg(F^{1,0}) =\\
g_0 \cdot \deg(\Omega^1_{Y}(\log S))=
{\rm rank}(F_1^{1,0})\cdot {\rm rank}(\sT)\cdot \deg(\Omega^1_{Y}(\log S)).
\end{gather*}
So $(F^{1,0}_1\oplus F^{0,1}_1,\theta_1)$ again satisfies the
assumptions made in \ref{semistableii}.

\begin{proof}[Proof of \ref{semistableii}]
Taking the determinant of
$$
\theta^{1,0}: F^{1,0}\> \simeq >> F^{0,1}\otimes\Omega^1_Y(\log
S),
$$
one obtains an isomorphism
$$
\det\theta^{1,0}: \det F^{1,0}\> \simeq >> \det
F^{0,1}\otimes{\Omega^1_Y(\log S)}^{g_0},
$$
By assumption there exists an invertible sheaf $\sL$ with
${\sL}^2=\Omega^1_Y(\log S)$. Since $F^{1,0}\simeq F^{0,1\vee},$
$$
(\det F^{1,0})^{2}\simeq {\Omega^1_Y(\log S)}^{g_0}={\sL}^{2\cdot g_0},
$$
and $\det F^{1,0}\otimes {\sL}^{-g_0}=\sN$ is of order two in
${\rm Pic}^0(Y)$.

If $g_0$ is even, $\sL$ is uniquely determined up to the tensor product with
two torsion points in ${\rm Pic}^0(Y)$.

If $g_0$ is odd, one replaces $\sL$ by $\sL\otimes \sN$
and obtains $\det F^{1,0}=\sL^{g_0}$.

By \ref{semistablei} the sheaf $$\sT=F^{1,0}\otimes \sL^{-1}$$ is
poly-stable of degree zero. \ref{simpson1} implies that the Higgs
bundle $(\sT,0)$ corresponds to a local system $\T$, necessarily
unitary.

Choose $\L$ to be the local system corresponding to the Higgs
bundle
$$(\sL\oplus\sL^{-1},\tau), \mbox{ \ \ with \ \ }
\tau^{1,0}: \sL \> \simeq >> \sL^{-1}\otimes \Omega^1_Y(\log S).
$$
The isomorphism
$$
\theta^{1,0}:\sT\otimes \sL = F^{1,0}\> \simeq >> F^{0,1}\otimes
\Omega^1_Y(\log S)\> \simeq >> F^{0,1}\otimes \sL^{2}
$$
induces an isomorphism
$$
\phi: \sT\otimes \sL^{-1}\> \simeq >> F^{0,1},
$$
such that $\theta^{1,0}=\phi\circ({\rm id}_\sT \otimes \tau^{1,0})$. Hence
the Higgs bundles $(F^{1,0}\oplus F^{0,1},\theta)$ and
$(\sT\otimes(\sL\oplus \sL^{-1}),{\rm id}_\sT\otimes\tau)$ are
isomorphic, and $\V\simeq \T\otimes_\C \L$.
\end{proof}

\begin{remark}\label{even} \
\begin{enumerate}
\item[i.] If $\deg \Omega^1_Y(\log S)$ is odd, hence $S\neq \emptyset$, and
if the genus of $Y$ is not zero, one may replace $Y$ by an \'etale
covering, in order to be able to apply \ref{semistableii}. Doing
so one may also assume that the invertible sheaf $\sN$ in
\ref{semistableii}, c), is trivial.
\item[ii.]
For $Y=\P^1$ and for $\X$ reaching the Arakelov bound, $\#S$ is
always even. This, together with the decomposition
\ref{semistableii}, for $\U=\C^{g_0}$, can easily obtained in the
following way. By \ref{semistablei}, $F^{1,0}$ must be the direct
sum of invertible sheaves $\sL_i$, all of the same degree, say
$\nu$. Since $\theta^{1,0}$ is an isomorphism, the image
$\theta^{1,0}(\sL_i)$ is $\sO_{\P^1}(2-s+\nu)\otimes \Omega$.
Since $F^{0,1}$ is dual to $F^{1,0}$ one obtains $-\nu=2-s+\nu$,
and writing $\sL_i^{-1}=\theta^{1,0}(\sL_i)$,
$$
(F^{1,0}\oplus F^{0,1}, \theta)\simeq (\bigoplus_i
\sO_{\P^1}(\nu)\oplus \sO_{\P^1}(-\nu), \bigoplus_i \tau).
$$
\end{enumerate}
\end{remark}

Consider now the local system of endomorphism $\END(\V)$ of $\V$,
which is a polarized weight zero variation of $L$ Hodge
structures. The Higgs bundle
$$
(F^{1,0}\oplus F^{0,1},\theta)
$$
for $\V$ induces the Higgs bundle
$$
(F^{1,-1}\oplus F^{0,0}\oplus F^{-1,1},\theta)
$$
corresponding to $\END(\V)=\V\otimes_\C\V^\vee$, by choosing
\begin{gather*}
F^{1,-1}=F^{1,0}\otimes {F^{0,1}}^\vee,
\ \ \ F^{0,0}=
F^{1,0}\otimes {F^{1,0}}^\vee \oplus F^{0,1}\otimes {F^{0,1}}^\vee\\
\mbox{ \ \ \ and \ \ \ }
F^{-1,1}=F^{0,1}\otimes {F^{1,0}}^\vee.
\end{gather*}
The Higgs field is given by
$$
\theta_{1,-1}= (-\id)\otimes{\tau_{1,0}}^\vee \oplus
\tau_{1,0}\otimes \id \mbox{ \ \ \ and \ \ \ }
\theta_{0,0}=\tau_{1,0}\otimes \id \oplus (-\id) \otimes
{\tau_{1,0}}^\vee.
$$
\begin{lemma}\label{split3}
Assume as in \ref{assumption} that $\X$ reaches the Arakelov bound
or equivalently that the Higgs field of $\V$ is maximal. Let
$$
F^{0,0}_u := {\rm Ker} (\tau_{0,0}) \mbox{ \ \ \ and \ \ \ }
F^{0,0}_m = {\rm Im}(\tau_{1,-1}).
$$
Then there is a splitting of the Higgs bundle
$$
(F^{1,-1}\oplus F^{0,0}\oplus F^{-1,1},\theta)=
(F^{1,-1}\oplus F_m^{0,0} \oplus F^{-1,1},\theta)\oplus (F^{0,0}_u,0),
$$
which corresponds to a splitting of the local system over $\C$
$$
\END(\V)=\W\oplus\U.
$$
$\U$ is unitary of rank $g_0^2$ and a variation of Hodge
structures concentrated in bidegree $0,0$, whereas $\W$ is a $\C$
variation of Hodge structures of weight zero and rank $3g_0^2$.
$$
\tau_{1,-1}:F^{1,-1}\>>> F^{0,0}_m\otimes \Omega^1_Y(\log S)
\mbox{ \ \ \ and \ \ \ } \tau_{0,0}:F^{0,0}_m\>>> F^{-1,1}\otimes
\Omega^1_Y(\log S)
$$
are both  isomorphisms.
\end{lemma}
\begin{proof}
By definition, $(F^{0,0}_u,0)$ is a sub Higgs bundle of
$(F^{1,-1}\oplus F^{0,0}\oplus F^{-1,1},\theta)$. We have an exact
sequence
$$
0\to F^{0,0}_u \to F^{0,0} \to F^{-1,1}\otimes \Omega_Y^1(\log S)
\to \sC
$$
where $\sC$ is a skyscraper sheaf. Hence
\begin{equation}\label{equat1.7}
\deg(F^{0,0}_u) \geq \deg(F^{0,0}) - \deg(F^{-1,1}) -
\rank(F^{-1,1})\cdot \deg(\Omega^1_Y(\log S)).
\end{equation}
Note that if (\ref{equat1.7}) is an equality then $\sC$ is
necessarily zero.

Since $\deg(F^{0,0})=0$ and since, by the Arakelov equality,
\begin{gather*}
\deg(F^{-1,1}) = g_0\cdot \deg(F^{0,1}) + g_0 \cdot
\deg({F^{1,0}}^\vee) \\
= g_0^2 \cdot \deg(\Omega_Y^1(\log S))= \rank(F^{-1,1}) \cdot
\deg(\Omega_Y^1(\log S))
\end{gather*}
one finds $\deg(F^{0,0}_u) \geq 0$. By \ref{simpson1} the degree
of $F_u^{0,0}$ can not be strictly positive, hence it is zero and
(\ref{equat1.7}) is an equality.

Again by \ref{simpson1} $(F^{0,0}_u,0)$ being a Higgs subbundle
of degree zero with trivial Higgs field, it corresponds to a
unitary local subsystem $\U$ of ${\END}(\V)$. The exact sequence
$$
0\to F^{0,0}_u\to F^{0,0}\to F^{-1,1}\otimes\Omega^1_Y(\log S)\to 0
$$
splits, and one obtains a direct sum decomposition of Higgs
bundles
$$
(F^{1,-1}\oplus F^{0,0}\oplus F^{-1,1},\theta)= (F^{1,-1}\oplus
F_m^{0,0} \oplus F^{-1,1},\theta)\oplus (F^{0,0}_u,0),
$$
which induces the splitting on $\END(\V)$ with the desired
properties.
\end{proof}
In \ref{split3} the local subsystem $\W$ of $\END(\V)$ has a
maximal Higgs field in the following sense.
\begin{definition}\label{maxhiggs} \
Let $\W$ be a $\C$ variation of Hodge structures of weight $k$,
and let
$$(F,\tau)=(\bigoplus_{p+q=k} F^{p,q},\bigoplus \tau_{p,q})$$
be the corresponding Higgs bundle. Recall that the width is
defined as
$$
{\rm width}(\W)={\rm Max}\{|p-q|; \ F^{p,q}\neq 0 \}.
$$
\begin{enumerate}
\item[i.]
$\W$ (or $(F,\tau)$) has a generically maximal Higgs field, if
${\rm width}(\W)>0$ and if
\begin{enumerate}
\item[a.] $F^{p,k-p} \neq 0$ for all $p$ with $|2p-k| \leq {\rm width}(\W)$.
\item[b.] $ \tau_{p,k-p}:F^{p,k-p} \to F^{p-1,k-p+1}\otimes\Omega^1_Y(\log S)$
is generically an isomorphism for all $p$ with $|2p-k| \leq {\rm
width}(\W)$ and $|2p-2-k| \leq {\rm width}(\W)$.
\end{enumerate}
\item[ii.] $\W$ (or $(F,\tau)$) has a maximal Higgs field, if the
$\tau_{p,k-p}$ in i), b. are all isomorphisms.
\end{enumerate}
\end{definition}
In particular, a variation of Hodge structures with a maximal Higgs field
can not be unitary.

\begin{properties}\label{width} \
\begin{enumerate}
\item[a.] If $\W$ is a $\C$ variation of Hodge structures with a
(generically) maximal Higgs field, and if
$\W'\subset \W$ is a direct factor, then ${\rm width}(\W')={\rm
width}(\W)$ and $\W'$ has again a (generically) maximal Higgs
field.
\item[b.] Let $\L$ and $\T$ be two variations of Hodge structures
with $\L\otimes\T$ of weight $1$ and width $1$, and with a
(generically) maximal Higgs field.

Then, choosing the bidegrees for $\L$ and $\T$ in an appropriate way,
either $\L$ is a variation of Hodge structures
concentrated in degree $0,0$, and $\T$ is a variation of Hodge structures
of weight one and width one with a (generically) maximal Higgs field,
or vice versa.
\end{enumerate}
\end{properties}
\begin{proof}
For a) consider the Higgs field $(\bigoplus{F'}^{p,q},
{\tau'}_{p,q}$ of $\W'$, which is a direct factor of the one for
$\W$. Since the $\tau_{p,q}$ are (generically) isomorphisms, a)
is obvious.

In b) denote the components of the Higgs fields of $\L$ and $\T$
by $\sL^{p_1,q_1}$ and $\sT^{p_2,q_2}$, respectively.
Shifting the bigrading one may assume that $p_1=0$ and $p_2=0$
are the smallest numbers with $\sL^{p_1,q_1}\neq 0$ and
$\sT^{p_2,q_2}\neq 0$ and moreover that the corresponding $q_i\geq 0$.
Since $q_1+q_2=1$, one of $q_i$ must be zero, let us say the first one.

Then $\sT^{p_2,q_2}$ can only be non-zero, for $(p_2,q_2)=(0,1)$
or $=(1,0)$ and $\L$ is concentrated in degree $0,0$.

Obviously this forces the Higgs field of $\L$ to be zero. Then
the Higgs field of $\L\otimes\T$ is the tensor product of the Higgs
field
$$
\sT^{1,0}\to \sT^{0,1}\otimes \Omega_Y^1(\log S)
$$
with the identity on $\sL^{0,0}$, hence the first one has to be
(generically) an isomorphism.
\end{proof}
\begin{remark}\label{altsplit}
The splitting in \ref{split3} can also be described by the tensor
product decomposition $\V= \T\otimes_\C\L$ in \ref{semistableii}
with $\T$ unitary and $\L$ a rank two variation of Hodge
structures of weight one and with a maximal Higgs field. For any
local system $\M$
one has a natural decomposition $\END(\M)=\END_0(\M) \oplus \C$,
where $\C$ acts on $\M$ by multiplication. Applying \ref{split3}
to $\L$ instead of $\V$, gives exactly the decomposition
$\END(\L)= \END_0(\L) \oplus \C$. One obtains
\begin{gather*}
\END(\V)= \T\otimes_\C {\T}^\vee \otimes_\C \L \otimes_\C {\L}^\vee =\\
(\END_0(\T)\oplus \C) \otimes_\C (\END_0(\L)\oplus \C)= \END_0(\T)
\oplus \C \oplus \END(\T)\otimes_\C \END_0(\L).
\end{gather*}
Here $\END_0(\T) \oplus \C$ is unitary and $\W=\END(\T)\otimes_\C
\END_0(\L)$ has again a maximal Higgs field.
\end{remark}

\begin{remark}\label{wedge}
If one replaces ${\END}(\V)$ by the isomorphic local system $\V
\otimes_\C \V$, one obtains the same decomposition. However, it
is more natural to shift the weights by two, and to consider this
as a variation of Hodge structures of weight $2$.

A statement similar to \ref{split3} holds true for
$\wedge^2(\V)$. Here the Higgs bundle is given by
\begin{gather*}
{F'}^{2,0}=F^{1,0}\wedge F^{1,0},
\ \ \ F^{1,1}=
{F'}^{1,0}\otimes F^{0,1}
\mbox{ \ \ \ and \ \ \ }
{F'}^{0,2}=F^{0,1}\wedge F^{0,1}.
\end{gather*}
\end{remark}
\section{Shimura curves and the special Mumford-Tate group}
\label{sectshimura}

\begin{lemma} \label{quotient}
Let $\L$ be a real variation of Hodge structures of weight $1$,
and of dimension $2$, with a non trivial Higgs field. Let
$\gamma_\L:\pi_1(U,*) \to \Sl(2,\R)$ be the corresponding
representation and let $\Gamma_\L$ denote the image of
$\gamma_\L$. Assume that the local monodromies around the points
$s\in S$ are unipotent. Then the Higgs field of $\L$ is maximal if and only
if $U=Y\setminus S\simeq \sH/\Gamma_\L$.
\end{lemma}
\begin{proof}

Writing $\sL$ for the $(1,0)$ part, we have an non trivial map
\begin{equation} \label{modular1}
\tau_{1,0}: \sL \>>> \sL^{-1} \otimes \Omega^1_{Y}(\log S).
\end{equation}
Since $\sL$ is ample, $\Omega^1_{Y}(\log S)$ is ample, hence the
universal covering $\tilde U$ of $U=Y\setminus S$ is the upper
half plane $\sH$. Let
$$
\tilde{U}= \sH \> \tilde \varphi >> \sH
$$
be the period map. The tangent sheaf of the period domain $\sH$
is given by the sheaf of homomorphisms from the $(1,0)$ part to
the $(0,1)$ part of the variation of Hodge structures. Therefore
$\tau_{1,0}$ is an isomorphism if and only if $\tilde \varphi$ is
a local diffeomorphism. Note that by Schmid \cite{Sch} the Hodge
metric on the Higgs bundle corresponding to $\L$ has logarithmic
growth at $S$ and bounded curvature. By the remarks following
\cite{Si3}, Propositions 9.8 and 9.1, $\tau_{1,0}$ is an
isomorphism if and only if $\tilde \varphi: \tilde U \to \sH$ is
a covering map, hence an isomorphism.

Obviously the latter holds true in case $Y\setminus S\simeq
\sH/\Gamma_\L.$

Assume that $\tilde\varphi$ is an isomorphism. Since $\tilde
\varphi$ is an equivariant with respect to the
$\pi_1(U,*)-$action on $\tilde U$ and the
$\text{P}\rho_{\L}(\pi_1(U,*))-$action on $\sH,$ the homomorphism
$$\rho_{\L_\Z}: \pi_1(U,*)\>>> \text{P}\rho_{\L_\Z}(\pi_1(U,*))\subset
\text{PSl}_2(\R)$$ must be injective, hence an isomorphism. So
$\tilde \varphi$ descend to an isomorphism
$$ \varphi: Y\setminus S\simeq \sH/\Gamma_\L.$$
\end{proof}

\begin{proof}[Proof of Proposition \ref{modular}]
$h:E\to Y$ be the semi-stable family of elliptic curves, reaching
the Arakelov bound, smooth over  $U$. Hence
$\L_{\mathbb Z}=R^1h_*\mathbb Z_{E_0}$ is a $\Z$-variation of
Hodge structures of weight one and of rank two. Writing $\sL$ for
the $(1,0)$ part, we have an isomorphism
\begin{equation} \label{modular2}
\tau_{1,0}: \sL \>>> \sL^{-1} \otimes \Omega^1_{Y}(\log S).
\end{equation}
Since $\sL$ is ample, $\Omega^1_{Y}(\log S)$ is ample, hence the
universal covering of $U$ is the upper half plane $\sH$. One
obtains a commutative diagram
$$
\begin{CD}
\sH \> \tilde \varphi >> \sH\\
\V \psi' VV \V \psi VV\\
U \> j >> \C
\end{CD}
$$
where $j$ is given by the $j$-invariant of the fibres of $E_0 \to
U$, where $\psi$ is the quotient map $\sH \to
\sH/\text{Sl}_2(\Z)$, and where $\tilde\varphi$ is the period map.
\ref{quotient} implies that
$$
\varphi:  U\>>> \sH/\rho_{\L_\Z}(\pi_1(U,*))
$$
is an isomorphism, hence $ \rho_{\L_\Z}(\pi_1(U,*))\subset
\text{Sl}_2(\Z)$ is of finite index, and  $E\to Y$ is a
semi-stable model of a modular curve.
\end{proof}

Let us recall the description of wedge products of tensor
products (see \cite{F-H}, p. 80). We will write
$\lambda=\{\lambda_1, \ldots , \lambda_\nu\}$ for the partition
of $g_0$ as $g_0=\lambda_1 + \cdots + \lambda_\nu$. The partition
$\lambda$ defines a Young diagram and the Schur functor ${\mathbb
S}_\lambda$. Assuming as in \ref{semistableii} that $\L$ is a
local system of rank $2$, and $\T$ a local system of rank $g_0$,
both with trivial determinant, one has
$$
\wedge^{k}(\L\otimes \T) = \bigoplus {\mathbb S}_\lambda (\L)
\otimes {\mathbb S}_{\lambda'}(\T)
$$
where the sum is taken over all partitions $\lambda$ of $k$ with
at most $2$ rows and at most $g_0$ columns, and where $\lambda'$
is the partition conjugate to $\lambda$. Similarly,
$$
S^{k}(\L\otimes \T) = \bigoplus {\mathbb S}_\lambda (\L) \otimes
{\mathbb S}_{\lambda}(\T)
$$
where the sum is taken over all partitions $\lambda$ of $k$ with
at most $2$ rows.

The only possible $\lambda$ are of the form $\{k-a,a\}$, for
$a\leq\frac{k}{2}$. By \cite{F-H}, 6.9 on p. 79,
$$
{\mathbb S}_{\{k-a,a\}}(\L)=\left\{ \begin{array}{lll} {\mathbb
S}_{\{k-2a\}}(\L)=S^{k-2a}(\L)\otimes\det(\L)^a & \mbox{ if } & 2a < k\\
{\mathbb S}_{\{1,1\}}(\L)=\det(\L)^a & \mbox{ if } & 2a=k
\end{array} \right. .
$$
For $k=g_0$ one obtains:
\begin{lemma}\label{tensor} Assume that $\det(\L)=\C$ and $\det(\T)=\C$.
\begin{enumerate}
\item[a.] If $g_0$ is odd, then for some partitions $\lambda_c$,
$$
\bigwedge^{g_0}(\L\otimes \T) = \bigoplus_{c=0}^{\frac{g_0-1}{2}}
S^{2c+1}(\L)\otimes {\mathbb S}_{\lambda_{2c}}(\T).
$$
In particular, for $c=\frac{g_0-1}{2}$ one obtains
$$
S^{g_0}(\L)\otimes \bigwedge^{g_0}(\T)=S^{g_0}(\L)
$$
as a direct factor.
\item[b.] If $g_0$ is even, then for some partitions $\lambda_c$,
$$
\bigwedge^{g_0}(\L\otimes \T) = S^{g_0}(\L)\oplus {\mathbb
S}_{\{2,\ldots,2\}}(\T) \oplus \bigoplus_{c=1}^{\frac{g_0}{2}-1}
S^{2c}(\L)\otimes {\mathbb S}_{\lambda_{2c}}(\T).
$$
\end{enumerate}
\end{lemma}

\begin{lemma}\label{tensor1} Assume that $\L$ and $\T$ are variations
of Hodge structures, with $\L$ of weight one, width one and with a
maximal Higgs field, and with $\T$ pure of bidegree $0,0$.
\begin{enumerate}
\item[a.] If $k$ is odd,
$$
H^0(Y,\bigwedge^k(\L\otimes\T)=0
$$
\item[b.] If $k$ is even, say $k=2c$, then for some $\lambda_c$
$$
H^0(Y,\bigwedge^k(\L\otimes\T)=H^0(Y,\det(\L)^c \otimes {\mathbb
S}_{\lambda_c}(\T)).
$$
\item[c.] For $k=2$ one has in ii) $\lambda_1=\{2\}$, hence ${\mathbb
S}_{\lambda_1}(\T))=\bigwedge^2(\T)$.
\item[d.] $H^0(Y,S^2(\L\otimes\T))= H^0(Y,\det(\L)\otimes\bigwedge^2(\T))$.
\end{enumerate}
\end{lemma}
\begin{proof}
$S^\ell(\L)$ has a maximal Higgs field for $\ell >0$, whereas for
all partitions $\lambda'$ the variation of Hodge structures
${\mathbb S}_{\lambda'}(\T)$ is again pure of bidegree $0,0$. By
\ref{width}, a), $S^\ell(\L)\otimes{\mathbb S}_{\lambda'}(\T)$
has no global sections. Hence $\bigwedge^k(\L\otimes \T)$ can
only have global sections for $k$ even. In this case, the global
sections lie in
$$
\det(\L)^c \otimes {\mathbb S}_{\lambda_c}(\T),
$$
for some partition $\lambda_c$, and one obtains a) and b). For
$k=2$ one finds $\lambda_1=\{2\}$. For d) one just has the two
partitions $\{1,1\}$ and $\{2\}$. Again, the direct factor
$S^2(\L)\otimes S^2(\T)$, corresponding to the first one, has no
global section.
\end{proof}

Let us shortly recall Mumford's definition of the Hodge group,
or as one writes today, the special Mumford-Tate group (see \cite{Mum1}, \cite{Mum2},
and also \cite{De2} and \cite{Scho}). Let $B$ be an abelian variety and
$H^1(B,\Q)$ and $Q$ the polarization on $V$. The special
Mumford-Tate group $\Hg(B)$ is defined in \cite{Mum1} as the
smallest $\Q$ algebraic subgroup of ${\rm Sp}(H^1(B,\R),Q)$,
which contains the complex structure. Equivalently $\Hg(B)$ is
the largest $\Q$ algebraic subgroup of ${\rm Sp}(H^1(B,\Q),Q)$,
which leaves all Hodge cycles of $B\times \cdots \times B$
invariant, hence all elements
$$
\eta\in H^{2p}(B\times \cdots \times B,\Q)^{p,p}
=\big[\bigwedge^{2p}(H^1(B,\Q)\oplus \cdots \oplus
H^1(B,\Q))\big]^{p,p}.
$$
For a smooth family of abelian varieties $f:X_0 \to U$ with
$B=f^{-1}(y)$ for some $y\in U$, and for the corresponding $\Q$
variation of polarized Hodge structures $R^1f_*\Q_{X_0}$ consider
Hodge cycles $\eta$ on $B$ which remain Hodge cycles under
parallel transform. One defines the special Mumford-Tate group
$\Hg(R^1f_*\Q_{X_0})$ as the largest $\Q$ subgroup which leaves
all those Hodge cycles invariant (\cite{De2}, \S 7, or
\cite{Scho}, 2.2).
\begin{lemma}\label{hg} \
\begin{enumerate}
\item[a.] For all $y\in U$ the special Mumford-Tate group $\Hg(f^{-1}(y))$
is a subgroup of $\Hg(R^1f_*\Q_{X_0})$. For all $y$ in the complement $U'$ of
the union of countably many proper closed subsets it coincides with
$\Hg(R^1f_*\Q_{X_0})$.
\item[b.] Let $G^{\rm Mon}$ denote the smallest reductive $\Q$ subgroup
of ${\rm Sp}(H^1(B,\R),Q)$, containing the image $\Gamma$ of the monodromy
representation
$$
\gamma:\pi_0(U)\>>> {\rm Sp}(H^1(B,\R),Q).
$$
Then the connected component $G^{\rm Mon}_0$ of one in
$G^{\rm Mon}$ is a subgroup of $\Hg(R^1f_*\Q_{X_0})$.
\item[c.] If $f:X\to Y$ reaches the Arakelov bound, and if $R^1f_*\C_X$
has no unitary part, then $G^{\rm Mon}_0=\Hg(R^1f_*\Q_{X_0})$.
\end{enumerate}
\end{lemma}
\begin{proof} The first statement of a) has been verified in \cite{Scho}, 2.3.,
and the second in \cite{Mo} 1.2. As explained in \cite{De2}, \S 7, or \cite{Scho}, 2.4, the
Mumford-Tate group contains a subgroup of $\Gamma$ of finite index,
hence b) holds true. It is easy to see, that the same holds true for
the special Mumford-Tate group (called Hodge group in \cite{Mum1}) by
using the same argument.

Since the special Mumford-Tate group of an abelian variety is
reductive, a) implies that $\Hg(R^1f_*\Q_{X_0})$ is reductive.
So $G^{\rm Mon}_0 \subset \Hg(R^1f_*\Q_{X_0})$ is an
inclusion of reductive groups. The proof of 3.1, (c), in
\cite{De3} carries over to show that both groups are equal,
if they leave the same tensors
$$
\eta\in
=\big[\bigwedge^{2p}(H^1(B,\Q)\oplus \cdots \oplus
H^1(B,\Q))\big]
$$
invariant.

Let $\eta\in H^k(B,\Q)$ be invariant under
$\Gamma$, and let $\tilde\eta$ be the corresponding global section of
$$
\bigwedge^{k}(R^1f_*\Q_{X_0})=\bigwedge^k(\L\otimes \T).
$$
By \ref{tensor1}, i) and ii), one can only have global sections
for $k=2c$, and those lie in
$$
\det(\L)^c \otimes {\mathbb S}_{\lambda_c}(\T).
$$
In particular they are of pure bidegree $c,c$.

The same argument holds true, if one replaces $B$ and $f:X\to Y$ by
any product, which implies c).
\end{proof}

For the Hodge group  $\Hg(R^1f_*\Q_X)=\Hg\subset {\rm Sp}(2g,\Q)$, as in Lemma \ref{hg}
Mumford considers the moduli functor $\mathcal M(\Hg)$ of isomorphy classes
of polarized abelian varieties with special Mumford-Tate group equal to
a subgroup of $\Hg$. He shows that $\mathcal M(\Hg)$ admits a quasi-projective
coarse moduli space $M(\Hg)$, which lies in the coarse moduli space of polarized
abelian varieties $A_g.$  By Mumford (\cite{Mum1}, Section 3, \cite{Mum2}, Sections 1-2)
$$ M(\Hg)=\Gamma \backslash\Hg(\R)/K $$
where $K$ is a maximal compact subgroup of $\Hg(\R),$ and $\Gamma$ an arithmetic subgroup
of $\Hg(\Q).$ The embedding $M(\Hg)\hookrightarrow A_g$ is a totally geodesic embedding,
and $M(\Hg)$ is a Shimura variety of Hodge Type $\Hg.$

Let $f: X_0\to U$ be a family of abelian varieties with the special Mumford-Tate group
$\Hg(R^1f_*\Q_{X_0})=\Hg.$  By Lemma \ref{hg}, a), $f$ induces a morphism
$$  U\to M(\Hg).$$
\begin{proof}[Proof of \ref{MTfamilies}]
By  Proposition \ref{semistableii} the image of the monodromy representation of $f$ lies in
$\Sl_2(\R)\times {\bf G}$, for some compact group $\bf G$,
and its Zariski closure is $Sl_2(\R)\times {\bf G}$. Hence,
$G_0^{\rm Mon}(\R)$ is again the product of $Sl_2(\R)$ with a compact group.
Lemma \ref{hg}, c), implies that
$$Hg(\R)=Sl_2(\R)\times {\bf G}'$$
for a compact group ${\bf G}'$, hence
$\Hg(\R)/K \simeq Sl_2(\R)/SO_2$ is the upper half plane $\sH$.

In particular, $\dim M(\Hg)=1$. Since we assumed the family to be non-isotrivial
and semi-stable, the morphism $U\to M(\Hg)$ is surjective.

Consider the composition $ \phi: U\to M(\Hg)\to A_g$.
Replacing $U$ by an \'etale covering, we may assume
that $X_0\to U$ is the pullback of a universal family of abelian varieties, defined
over an \'etale covering $A'_g$ of $A_g$. The pull back of the tangent bundle on
$A_g$ via $\phi$ is just
$$\phi^*T_{A_g}=S^2E^{0,1}\subset {E^{0,1}}^{\otimes 2}.$$
The differential $d\phi: T_{U}\to \phi^*T_{A_g}\subset {E^{0,1}}^{\otimes 2}$ is induced
by the Kodaira-Spencer map $ E^{1,0}\otimes T_{U}\to E^{0,1}.$
By Proposition \ref{semistableii}
$$ E^{1,0}\oplus E^{0,1}=(\sL\oplus\sL^{-1})\otimes\T,$$
and the map $d\phi: T_{U}\to {E^{0,1}}^{\otimes 2}$ lies in the component
$$ d\phi: T_{U}\simeq\sL^{\otimes-2}\subset \sL^{\otimes -2}\otimes \rm{End}(\T).$$
This implies  that the differential of the map $U\to M(\Hg)$ is no where vanishing,
hence $U\to M(\Hg)$ is \'etale.
\end{proof}

\begin{remarks}\label{remark} \ 
\begin{enumerate}
\item[a)] As well known (see \cite{Mum1}, \cite{Mum2}) the moduli space of abelian varieties with a given
special Mumford-Tate group is necessarily a Satake holomorphic embedding. Hence the assumptions
made in Proposition \ref{MTfamilies} imply in particular that the period map from $U$ to the
corresponding moduli space of abelian varieties with a fixed level structure is
a Satake holomorphic embedding. 
\item[b)] Presumably Proposition \ref{MTfamilies} can also be obtained using 
\cite{Abd}. Using Proposition \ref{semistableii} the maximality of the Higgs field 
should imply that the period map from $U=\sH/\Gamma$ to the Siegel upper half plane is a rigid, 
totally geodesic, and equivariant holomorphic map. Then \cite{Abd}, Theorem 3.4, 
implies that $f: X\to Y$ is a family of Mumford type, and as mentioned in the introduction
one can finish the proof of Theorem \ref{shimura} going through the classification
of Shimura varieties.
\item[c)] Without the assumption of rigidity, hidden behind the one saying that the maximal
unitary local subsystem is defined over $\Q$, we do not see a way to show directly, that
the families are Kuga fibre spaces. One needs a precise description of the
$\Z$ structure on the decompositions of the variation of Hodge structures.
On the other hand, the latter will allow to prove Theorems \ref{shimura} and \ref{shimura2}
directly.
\item[d)] Theorems \ref{shimura} and \ref{shimura2} imply
that all families $f: X\to Y$ with maximal Higgs fields are Kuga fibre spaces, and that 
the period map is again a Satake holomorphic embedding.
\end{enumerate}
\end{remarks}

\section{Splitting over $\bar{\Q}$}
\label{sectbarsplitting}

Up to now, we considered local systems of $\C$-vector
spaces induced by the family of abelian varieties. We
say that a $\C$ local system $\M$ is defined over a subring $R$
of $\C$, if there exists a local system $\M_R$ of torsion free $R$-modules
with $\M=\M_R\otimes_R\C$. In different terms, the representation
$$
\gamma_\M:\pi_0(U,*) \>>> {\rm Gl}(\mu,\C)
$$
is conjugate to one factoring like
$$
\gamma_\M:\pi_0(U,*) \>>> {\rm Gl}(\mu,R) \>>> {\rm Gl}(\mu,\C).
$$
If $\M$ is defined over $R$, and if $\sigma:R\to R'$ is an
automorphism, we will write $\M_R^\sigma$ for the local system
defined by
$$
\gamma_\M:\pi_0(U,*) \>>> {\rm Gl}(\mu,R) \> \sigma >> {\rm
Gl}(\mu,R'),
$$
and $\M^\sigma=\M_R^\sigma\otimes_{R'} \C$. In this section we
want to show, that the splittings $\X=\V\oplus \U_1$ and
$\END(\V)=\W\oplus\U$ considered in the last section are defined
over $\bar{\Q}$, i.e. that there exists a number field $K$
containing the field of definition for $\X$ and local $K$ subsystems
$$
\V_K\subset \X_K, \ \U_{1K}\subset\X_K, \ \W_K\subset \END(\X_K)
\mbox{ \ and \ } \U_K\subset\END(\X_K)
$$
with
\begin{gather*}
\X_K = \X_L\otimes_LK = \V_K\oplus \U_{1K}, \ \ \V_K=\W_K \oplus \U_K,
\mbox{ \ \ and with}\\
\V=\V_K\otimes_K \C, \ \ \U_1=\U_{1K}\otimes_K\C \ \ \W=\W_K\otimes_K \C, \ \
\U=\U_K\otimes_K \C.
\end{gather*}

We start with a simple observation. Suppose that $\M$ is a local
system defined over a number field $L$. The local system $\M_{L}$
is given by a representation $\rho: \pi_1(U,*) \to
\text{Gl}(M_{L})$ for the fibre $M_{L}$ of $\M_{L}$ over the base
point $*$.

Fixing a positive integer $r,$ let $\mathcal G(r,\M)$ denote the
set of all rank-r local subsystems of $\M$ and let ${\rm
Grass}(r,M_L)$ be the Grassmann variety of $r$-dimensional subspaces.
Then $\mathcal G(r,\M)$ is the subvariety of
$$
{\rm Grass}(r, M_{L})\times_{{\rm Spec}(L)}{\rm Spec}(\C)
$$
consisting of the $\pi_1(U,*)$ invariant points. In particular,
it is a projective variety defined over $L$. An $K$-valued point
of $\mathcal G(r,\M)$ corresponds to a local subsystem of
$\M_K=\M_{L}\otimes_L K$. One obtains the following well known property.

\begin{lemma}\label{fieldofdef}
If $[\W]\in \mathcal G(r,\M)$ is an isolated point,
then $\W$ is defined over $\bar\Q.$
\end{lemma}

In the proof of  \ref{split6} we will also need:

\begin{lemma}\label{fieldofdef2}
Let $\M$ be an variation of Hodge structures
defined over $L$, and let $\W\subset \M$ be an irreducible local subsystem
of rank $r$ defined over $\C,$. Then $\W$ can be deformed to a local subsystem
$\W_t\subset\M$, which is isomorphic to $\W$ and which is defined over
a finite extension of $L$.
\end{lemma}
\begin{proof} By \cite{Del} $\M$ is completely reducible over $\C$. Hence
we have a decomposition $\M=\W\oplus\W'$.

The space $\mathcal G(r,\M)$ of rank $r$ local subsystems of $\M$
is defined over $L$ and the subset
$$ \{\W_t\in \mathcal G(r,\M); \mbox{ the composit }
\W_t \subset \W\oplus \W' \> pr_1 >> \W\mbox{ is non zero }\}$$
forms a Zariski open subset. So there exists some $\W_t$ in this subset,
which is defined over some finite extension of $L$. Since $p: \W_t\to \W$
is non zero, ${\rm rank}(\W_t)= {\rm rank}(\W)$, and since
$\W$ is irreducible, $p$ is an isomorphism.
\end{proof}
\begin{lemma}\label{splitqbar}
Let $\M$ be the underlying local system of a variation of Hodge
structures of weight $k$ defined over a number field $L$.
Assume that there is a decomposition
\begin{equation}\label{split5b}
\M=\U\oplus \bigoplus_{i=1}^\ell \M_i
\end{equation}
in sub variations of Hodge structures, and let
\begin{equation}\label{split5a}
(E,\theta)=(N,0)\oplus \bigoplus_{i=1}^\ell
(F_i,\tau_i=\theta|_{F_i})
\end{equation}
be the induced decomposition of the Higgs field. Assume that ${\rm
width}(\M_i)=i$, and that the $\M_i$ have all generically maximal
Higgs fields. Then the decomposition (\ref{split5b}) is defined
over $\bar\Q$. If $L$ is real, it is defined over
$\bar{\Q}\cap \R$. If $\M$ is polarized, then the decomposition
(\ref{split5b}) can be chosen to be orthogonal with respect to
the polarization.
\end{lemma}

\begin{proof}
Consider a family $\{\W_t\}_{t\in \Delta}$ of local subsystems of
$\M$ defined over a disk $\Delta$ with $\W_0=\M_\ell$. For $t\in
\Delta$ let $(F_{\W_t},\tau_t)$ denote the Higgs bundle
corresponding to $\W_t$. Hence $(F_{\W_t},\tau_t)$ is obtained by
restricting the $F$-filtration of $\M\otimes \sO_U$ to
$\W_t\otimes\sO_U$ and by taking the corresponding graded sheaf.
So the Higgs map
$$ \tau^{p,k-p}: F^{p,k-p}_t\>>>
F^{p-1,k-p+1}_t \otimes\Omega^1_Y(\log S)$$ will again be
generically isomorphic for $t$ sufficiently closed to $0$ and
$$
|2p-k| \leq \ell \mbox{ \ \ and \ \ } |2p-2-k| \leq \ell.
$$
If the projection
$$\rho: \W_t \>>> \M = \U\oplus \bigoplus_{i=1}^\ell\M_i
\>>> \U\oplus \bigoplus_{i=1}^{\ell-1} \M_i$$ is non-zero, the
complete reducibility of local systems coming from variations of
Hodge structures (see \cite{Del}) implies that $\W_t$ contains an
irreducible non-trivial direct factor, say $\W'_t$ which is
isomorphic to a direct factor of $\U$ or of one of the local
systems $\M_i$, for $i <  \ell$.

Restricting again the $F$ filtration and passing to the
corresponding graded sheaf, we obtain a Higgs bundle
$(F_{\W'_t},\tau'_t)$ with trivial Higgs field, or whose width is
strictly smaller than $\ell$. On the other hand,
$(F_{\W'_t},\tau'_t)$ is a sub Higgs bundle of the Higgs bundle
$(F_{\W_t},\tau_t)$ of width $\ell$, a contradiction. So $\rho$
is zero and $\W_t=\M_{\ell}$.

Thus $\M_{\ell}$ is rigid as a local subsystem of $\M$, and by
Lemma \ref{fieldofdef} $\M_{\ell}$ is defined over $\bar{\Q}$.

Assume now that $L$ is real, hence $\M=\M_\R\otimes \C$.
The complex conjugation
defines an involution $\iota$ on $\M$. Let $\bar\M_{\ell}$ denote
the image of $\M_{\ell}$ under $\iota$. Then $\bar\M_{\ell}$ has
again generically isomorphic Higgs maps $\tau^{p,k-p},$ for
$$
|2p-k| \leq \ell \mbox{ \ \ and \ \ } |2p-2-k| \leq \ell.
$$
If $\bar\M_{\ell}\neq\M_{\ell}$, repeating the argument used
above, one obtains a map
$$
\bar \M_{\ell}\>>> \U\oplus \bigoplus_{i=1}^{\ell-1} \M_i,
$$
from a Higgs bundle of width $\ell$ and with a maximal Higgs field
to one with trivial Higgs field or of lower width. Again such a
morphism must be zero, hence $\M_\ell = \bar\M_\ell$ in this case.

So we can find a number field $K$, real in case $L$ is real,
and a local system
$\M_{\ell,K} \subset \M_K$ with $\M_{\ell}=\M_{\ell,K}\otimes_K
\C$. The polarization on $\M$ restricts to a non-degenerated
intersection form on $\M_K$. Choosing for $\M_{\ell,K}^{\perp}$
the orthogonal complement of $\M_{\ell,K}$ in $\M_K$ we obtain a
splitting
$$ \M_K=\M_{\ell,K}\oplus \M_{\ell,K}^\perp$$
inducing over $\C$ the splitting of the factor $\M_{\ell}$ in
(\ref{split5b}). By induction on $\ell$ we obtain \ref{splitqbar}.
\end{proof}

For a reductive algebraic group $G$ and for a finitely generated
group $\Gamma$ let $\sM(\Gamma, G)$ denote the moduli space of
reductive representations of $\Gamma$ in $G$.

\begin{theorem}[Simpson, \cite{Si4}, Cor.9.18]
\label{simpson2} \ Suppose $\Gamma$ is a finitely generated
group. Suppose $\phi: G\to H$ is a homomorphism of reductive
algebraic groups with finite kernel. Then the resulting morphism
of moduli  spaces
$$\phi: \mathcal M(\Gamma, G)\>>> \mathcal M(\Gamma, H)$$
is finite.
\end{theorem}

\begin{corollary}\label{simpson3} Let $\Gamma$ be $\pi_1(Y,*)$
of a projective manifold, and $\gamma:\Gamma\to G$ be a reductive
representation. If $\phi\gamma\in \mathcal M(\Gamma, H)$ comes
from a $\C$ variation of Hodge structures, then $\gamma$ comes from a
$\C$ variation of Hodge structures as well.
\end{corollary}

\begin{proof} By Simpson a reductive local system is coming from an
variation of Hodge structures if and only if the isomorphism class
of the corresponding Higgs bundle is a fix point of the $\C^*$
action. Since the $\C^*$ action contains the identity and since it
is compatible with $\phi$, the finiteness of the preimage
$\phi^{-1}\phi(\gamma)$ implies that the isomorphism class of the
Higgs bundle corresponding to  $\gamma$ is fixed by the $\C^*$
action, as well.
\end{proof}

\begin{definition}\label{deftrace}
Let $\M$ be a local system of rank $r$, and defined over $\bar{\Q}$.
Let $\gamma_\M:\pi_1(U,*) \to \Sl(2, \bar{\Q})$ be the corresponding
representation of the fundamental group. For $\eta\in \pi_1(U,*)$ we write
$\tr(\gamma_\M(\eta))\in \bar{\Q}$ for the trace of $\eta$ and
$$
\tr(\M)=\{\tr(\gamma_\M(\eta)); \ \eta \in \pi_1(U,*)\}.
$$
\end{definition}

\begin{corollary}\label{split6}
Under the assumptions made in \ref{assumption}
\begin{enumerate}
\item[i.] The splitting $\X=\V\oplus\U_1$ is defined over $\bar\Q$,
and over $\bar\Q \cap \R$ in case $L$ is real.
If $\X$ is polarized, it can be chosen to be orthogonal.
\item[ii.] The splitting $\END(\V)=\W\oplus\U$ constructed in
Lemma \ref{split3} is defined over $\bar\Q$, and over $\bar\Q\cap\R$
in case $L$ is real. If $\X$
is polarized, it can be chosen to be orthogonal.
\item[iii.] Replacing $Y$ by an \'etale covering $Y'$,
one can choose the decomposition $\V\simeq \L\otimes\T$ in
\ref{semistableii} such that
\begin{enumerate}
\item[a.] $\L$ and $\T$ are defined over a number field $K$, real
if $L$ is real.
\item[b.] One has an isomorphism $\V_{\bar\Q}\simeq
\L_{\bar\Q}\otimes_{\bar\Q}\T_{\bar\Q}$.
\item[c.] $\tr(\L)$ is a subset of the ring of integers $\sO_K$ of $K$.
\end{enumerate}
\end{enumerate}
\end{corollary}
\begin{proof}
i) and ii) are direct consequences of \ref{splitqbar}. For iii)
let us first remark that for $L$ real, passing to an \'etale
covering $\L$ and $\T$ can both be assumed to be defined over
$\R$. In fact, the local system $\bar{\L}$ has a maximal Higgs
field, hence its Higgs field is of the form $(\sL'\oplus
{\sL'}^{-1}, \tau')$ where $\sL'$ is a theta characteristic.
Hence it differs from $\sL$ at most by the tensor product with a
two torsion point in ${\rm Pic}^0(Y)$. Replacing $Y$ by an
\'etale covering, we may assume $\L=\bar{\L}$. From
\ref{semistableii}, d), we obtain $\bar{\T}=\T$.

Consider the isomorphism of local systems $\phi: \L \otimes \T \>
\simeq >> \V$ and the induced isomorphism
$$
\phi^2: \END_0(\L \otimes \T) = \END_0(\L) \oplus \END_0(\T)\otimes \END_0(\L)
\oplus \END_0(\T) \> >> \END(\V).
$$
Since $\phi^2 \END_0(\T)$ is the unitary part of this decomposition, by \ref{splitqbar}
it is defined over $\bar{\Q}\cap\R$, as well as
$\phi^2(\END_0(\L) \oplus \END_0(\T)\otimes \END_0(\L))$.
The $1,-1$ part of the Higgs field corresponding to $\phi^2\End_0(\L)$ has rank one,
and its Higgs field is maximal. Hence $\phi^2\End_0(\L)$ is irreducible, and by
\ref{fieldofdef2} it is isomorphic to a local system, defined over $\bar{\Q}$.
Hence $\T\otimes \T \simeq \END(\T)$ and $\L\otimes\L\simeq\END(\L)$ are
both isomorphic to local systems defined over some real number field $K'$.
An $\sO_{K'}$ structure can be defined by
$$
\phi^2(\END(\L))_{\sO_{K'}} = \phi^2(\END(\L))_{K} \cap \V_{\sO_{K'}}
$$
Consider for $\nu=2$ or $\nu=g_0$ the moduli
space $\sM (U, \Sl(\nu^2))$ of reductive representations of
$\pi(U,*)$ into $Sl(\nu^2)$.
It is a quasi-projective variety defined over $\Q$. The fact that
$\L\otimes\L$ (or $\T\otimes\T$) is defined over $\bar\Q$ implies that its
isomorphy class in $\sM(U, \Sl(\nu^2))$ is a $\bar\Q$ valued point.

Consider the morphism induced by the second tensor product
$$
\rho:\sM(U, \Sl(\nu)) \>>> \sM(U, \Sl(\nu^2))
$$
which is clearly defined over $\Q$. By \ref{simpson2}
$\rho$ is finite, hence the fibre $\rho^{-1}([\L\otimes\L])$ (or
$\rho^{-1}([\T\otimes \T])$) consists of finitely many
$\bar{\Q}$-valued points, hence $\L$ and $\T$ can be defined over
a number field $K$. If $L$ is real, as already remarked above,
we may choose $K$ to be real.

Obviously, for $\rho\in \pi_1(Y,*)$ one has
$$\tr(\gamma_\L(\rho))^2=\tr(\gamma_{\L\otimes \L}(\rho)).$$
In fact, one may assume that $\gamma_\L(\rho)$ is a diagonal
matrix with entries $a$ and $b$ on the diagonal. Then
$\tr(\gamma_{\L\otimes \L}(\rho))$ has $a^2$, $b^2$, $ab$ and $ba$
as diagonal elements. Since $\tr(\gamma_{\L\otimes \L}(\rho))\in
\sO_{K'}$ we find $\tr(\gamma_\L(\rho))\in \sO_K$.
\end{proof}

\section{Splitting over $\Q$ for $S\neq \emptyset$ and isogenies}
\label{sectqsplitting}

In this section, we will consider the case $L=\Q$ and $\X_\Q=R^1f_*\Q_{X_0}$,
where $f:X\to Y$ is a family of abelian varieties,
$S=Y\setminus U \neq \emptyset$, and where the restriction $X_0\to U$ of $f$ is a smooth family.
\begin{lemma}\label{split7}
Assume that $S\neq \emptyset$ and let $\M_\Q$ be a $\Q$-variation
of Hodge structures of weight $k$ and with unipotent monodromy
around all points $s\in S$. Assume that over some number field $K$
there exists a splitting
$$
\M_K=\M_\Q\otimes_\Q K=\W_K\oplus\U_K
$$
where $\U=\U_K\otimes_K\C$ is unitary and where the Higgs field of
$\W=\W_K\otimes_K\C$ is maximal. Then $\W$, $\U$ and the
decomposition $\M=\W\oplus\U$ are defined over $\Q$.
Moreover, $\U$ extends to a local system over $Y$.
\end{lemma}
\begin{proof}
Let $\T$ be a local subsystem of $\W$. Writing
$$
\big( \bigoplus_{p+q=k}F_{\T}^{p,q}, \bigoplus_{p+q=k}
\theta_{p,q}\big),
$$
for the Higgs bundle corresponding to $\T$, the maximality of the
Higgs field for $\W$ implies that the Higgs field for $\T$ is
maximal, as well. In particular, for all $s\in S$ and for $p>0$
the residue maps
$$
{\rm res}_s(\theta_{p,q}) : F_{\T,s}^{p,q} \>>> F_{\T,s}^{p-1,q+1}
$$
are isomorphisms. By \cite{Si2} the residues of the Higgs field at
$s$ are defined by the nilpotent part of the local monodromy
matrix around $s$. Hence if $\gamma$ is a small loop around $s$ in
$Y$, and if $\rho_\T(\gamma)$ denotes the image of $\gamma$ under
a representation of the fundamental group, defining $\T$, the
nilpotent part $N(\rho_\T(\gamma))=\log\rho_\T(\gamma)$ of $\rho_\T(\gamma)$
has to be non-trivial

We may assume that $K$ is a Galois extension of $\Q$. Recall that
for $\sigma \in {\rm Gal}(K/\Q)$ we denote the local systems
obtained by composing the representation with $\sigma$ by an
upper index $\sigma$. Consider the composite
$$
p:\U_K^\sigma \>>> \M_K=\W_K\oplus \U_K \>>> \W_K,
$$
and the induced map $\U^\sigma=\U_K^\sigma\otimes_K\C \to \W$.

Let $\gamma$ be a small loop around $s\in S$, and let
$\rho_\U(\gamma)$ and $\rho_{\U^\sigma}$ be the images of $\gamma$
under the representations defining $\U$ and $\U^\sigma$
respectively. Since $\U$ is unitary and unipotent, the nilpotent
part of the monodromy matrix $N(\rho_{\U}(\gamma))=0$. This being
invariant under conjugation, $N(\rho_{\U^\sigma}(\gamma))$ is
zero, as well as $N(\rho_{p(\U^\sigma)}(\gamma))$.

Therefore $p(\U^\sigma))=0$, hence $\U^\sigma=\U$, and $\U$ is
defined over $\Q$. Taking again the orthogonal complement, one
obtains the $\Q$-splitting asked for in \ref{split7}.

Since $N(\rho_{\U}(\gamma))=0$, the residues of $\U$ are zero in
all points $s\in S$, hence $\U$ extends to a local system on $Y$.
\end{proof}

\begin{corollary}\label{split8} Suppose that $S\neq\emptyset.$ Then
the splittings in Corollary \ref{split6}, i) and ii), can be defined over
$\Q.$
\end{corollary}

\begin{lemma}\label{unitary}
Let $\M$ be a local system, defined over $\Z$, and let
$\M_\Q=\W_\Q\oplus \U_\Q$ be a decomposition, defined over $\Q$.
Then there exist local systems $\U_\Z$ and $\W_\Z$, defined over
$\Z$ with
\begin{equation}\label{zstructure}
\U_\Q=\U_\Z\otimes \Q, \ \ \ \W_\Q=\W_\Z\otimes \Q, \ \ \ \mbox{
\ \ and \ \ } \M_\Z \supset \W_\Z \oplus \U_\Z.
\end{equation}
Moreover, if $\U_\Q$ is unitary with trivial local monodromies
around $S$, then there exists an \'etale covering $\pi:Y'\to Y$
such that $\pi^*\U_\Q$ is trivial.
\end{lemma}
\begin{proof}
Defining a $\Z$ structure on $\W_\Q$ and $\U_\Q$ by
$$
\W_\Z=\W_\Q\cap \M_\Z \mbox{  \ \ and \ \ } \U_\Z = \U_\Q \cap
\M_\Z
$$
(\ref{zstructure}) obviously holds true.

Since the integer elements of the unitary group form a finite
group, the representation defining $\U$ factors through a finite
quotient of the fundamental group $\pi_1(U,*)\mapsto G$. The
condition on the local monodromies implies that this quotient
factors through $\pi_1(Y,*)$, and we may choose $Y'$ to be the
corresponding \'etale covering.
\end{proof}
By \ref{split8} we obtain decompositions
$$
R^1f_*\Q_{X_0}=\V_\Q\oplus \U_{1\Q}\mbox{ \ \ and \ \ }
\END(\V_\Q)=\W_\Q \oplus \U_\Q.
$$
By \ref{split7} the local monodromies of the unitary parts $\U_1$
and $\U$ are trivial. Moreover, $\U$ is a sub variation of Hodge
structures of weight $0,0$. Summing up, we obtain:
\begin{corollary}\label{split9}
Let $f:X\to Y$ be a family of abelian varieties with unipotent local monodromies
around $s\in S$, and reaching the Arakelov bound. If $S\not=\emptyset$ there exists a
finite \'etale cover $\pi:Y'\to Y$  with
\begin{enumerate}
\item[i.]
$\pi^*(R^1f_*(\Z_{X_0}))\supset \V'_{\Z}\oplus \Z^{2(g-g_0)}$, and\\[.1cm]
$\pi^*(R^1f_*(\Z_{X_0}))\otimes\Q=(\V'_{\Z}\oplus \Z^{2(g-g_0)})
\otimes\Q,$\\[.1cm]
where $\V'_{\Z}$ is an $\Z$-variation of Hodge structures of
weight $1$ with maximal Higgs field.
\item[ii.]
$\END(\V'_{\Z})\supset \W'_{\Z}\oplus\Z^{g_0^2},\quad
\END(\V'_{\Z})\otimes\Q=(\W'_{\Z}\oplus\Z^{g_0^2})\otimes\Q,$\\[.1cm]
where $\W'_{\Z}$ is an $\Z$-variation of Hodge structures of
weight 0 with maximal Higgs field, and where $\Z^{g_0^2}$ is a
local $\Z$ subsystem of type $(0,0)$.
\end{enumerate}
\end{corollary}

\begin{proof}[Proof of Theorem \ref{geomsplit}]
Let $Y'$ be the \'etale covering constructed in \ref{split9},
ii). So using the notations introduced there,
\begin{equation}\label{zstructure2}
R^1f'_*(\Z_{X'_0})\otimes \Q = \V'_\Q \oplus \Z^{2(g-g_0)} \mbox{\
\ and \ \ } \END(\V'_\Q) = \W'_\Q \oplus \Z^{g_0^2}.
\end{equation}
The left hand side of (\ref{zstructure2}) implies that $f':X'\to
Y'$ is isogenous to a product of a family of $g_0$ dimensional
abelian varieties with a constant abelian variety $B$ of
dimension $g-g_0$. By abuse of notations we will assume from now
on, that $B$ is trivial, hence $g=g_0$ and
$R^1f'_*(\Z_{X'_0})\otimes \Q = \V'_\Q$, and we will show that
under this assumption $f':X'\to Y'$ is isogenous to a $g$-fold
product of a modular family of elliptic curves.

Let us write
$$\End(*)=H^0(Y',\END(*))
$$
for the global endomorphisms. $\End(\V'_\Q)=\Q^{g^2}$ is a $\Q$
Hodge structure of weight zero, in our case the Hodge filtration
is trivial, i.e.
$$
\End(\V'_\Q)^{0,0}=\End(\V'_\Q).
$$
If $X_\eta=X'\times_{Y'}{\rm Spec}(\overline{\C(Y')})$ denotes the general
fibre of $f'$, one obtains from \cite{Del}, 4.4.6,
$$
\End(X_\eta)\otimes \Q = \End(\V'_\Q)^{0,0}=\End(\V'_\Q).
$$

By the complete reducibility of abelian varieties,
there exists simple abelian varieties $B_1,\ldots,B_r$
of dimension $g_i$, respectively,
which are pairwise non isogenous, and such that
$X_\eta$ is isogenous to the product
$$
B_1^{\times \nu_1} \times \cdots \times B_r^{\times \nu_r}.
$$
Moreover, since $\V$ has no flat part, none of the $B_i$
can be defined over $\C$. Let us assume that
$g_i=1$ for $i=1,\ldots ,r'$ and $g_i>1$ for
$i=r'+1,\ldots ,r$.

By \cite{Mum}, p. 201, $D_i=\End(B_i)\otimes\Q$ is a division algebra
of finite rank over $\Q$ with center $K_i$.
Let us write
$$
d_i^2=\dim_{K_i}(D_i) \mbox{ \ \ and \ \ }
e_i=[K_i:\Q].
$$
Hence $e_i\cdot d_i^2=\dim_\Q(D_i)$.

By \cite{Mum}, p. 202, or by \cite{BL}, p. 141, either $d_i \leq
2$ and $e_i \cdot d_i$ divides $g_i$, or else $e_i\cdot d_i^2$
divides $2\cdot g_i$. In both cases, the rank $e_i\cdot d_i^2$ is
smaller than or equal to $2\cdot g_i$. If $i\leq r'$, hence if
$B_i$ is an elliptic curve, not defined over $\C$, we have
$e_i=d_i=1$.

Writing $M_{\nu_i}(D_i)$ for the $\nu_i\times \nu_i$ matrices
over $D_i$, one finds (\cite{Mum}, p. 174)
$$
\End(X_\eta)\otimes \Q =
M_{\nu_1}(D_1)\oplus\cdots \oplus M_{\nu_r}(D_r)
$$
hence
\begin{multline*}
g^2\leq \dim_\Q(\End(X_\eta)\otimes \Q)=
\big(\sum_{i=1}^r \nu_i\cdot g_i \big)^2 =
\sum_{i=1}^r (e_i\cdot d_i^2) \cdot \nu_i^2
\leq\\
\sum_{i=1}^{r'}\nu_i^2 + \sum_{i=r'+1}^r \nu_i^2\cdot 2\cdot g_i
\leq \sum_{i=1}^r \nu_i^2\cdot g_i^2 .
\end{multline*}
Obviously this implies that $r=1$ and that $g_1 \leq 2$. If
$g_1=1$, we are done. In fact, the isogeny extends all over
$Y'\setminus S'$ and, since we assumed the monodromies to be
unipotent, $B_1$ is the general fibre of a semi-stable family of
elliptic curves. The Higgs field for this family is again maximal,
and \ref{geomsplit} follows from \ref{modular}.\\

It remains to exclude the case that $g_1=2$, and that $e_1\cdot d_1^2=4$.
If the center $K_1$ is not a totally real number field, $e_1$ must be
lager than $1$ and one finds\\[.1cm]
I. $d_1=1$ and $D_1=K_1$ is a quadratic imaginary extension of a
real quadratic extension of $\Q$.\\[.1cm]
If $K_1$ is a real number field, looking again to the
classification of endomorphisms of simple abelian varieties in
\cite{Mum} or \cite{BL}, one finds that $e_1$ divides $g_1$, hence
the only possible case is \\[.1cm]
II. $d_1=2$ and $e_1=1$, and $D_1$ is a quaternion algebra over
$\Q$.\\[.1cm]
The abelian surface $B_1$ over ${\rm Spec}(\C(Y'))$ extends to a
non-isotrivial family of abelian varieties $B'\to Y'$, smooth outside of
$S$ and with unipotent monodromies for all $s\in S$. This family
again has a maximal Higgs field, and thereby the local monodromies
in $s \in S$ are non-trivial. As we will see below, in both cases,
I and II, the moduli scheme of abelian surfaces with the corresponding type of
endomorphisms turns out to be a compact subvariety of the
moduli scheme of polarized abelian varieties, a contradiction.\\[.1cm]
I. By \cite{BL}, Example 6.6 in Chapter 9, there are only finitely
many $g_1$ dimensional abelian varieties with a given type of
complex multiplication, i.e. with $D_1$ a quadratic imaginary
extension of a real number field of degree $g_1$ over $\Q$.\\[.1cm]
II. By \cite{BL}, Exercise (1) in Chapter 9, there is no abelian
surface for which $D_1$ is a totally definite quaternion algebra.
If $D_1=\End(B)\otimes\Q$ is totally indefinite, B is a false elliptic curve,
as considered in Example \ref{mumfordshimura}, ii, for $d=1$.
Such abelian surfaces have been studied in \cite{Sha},
and their moduli scheme is a compact Shimura curve.
The latter follows from Shimura's construction of the moduli
scheme as a quotient of the upper half plane $\sH$
(see \cite{BL}, \S 8 in Chapter 9, for example)
and from \cite{Shi}, Chapter 9.
\end{proof}

\section{Quaternion algebras and Fuchsian groups}
\label{sectquaternion}

Let $A$ denote a quaternion algebra over a totally real algebraic
number field $F$ with $d$ distinct embeddings
$$\sigma_1=id,\,\sigma_2,\ldots , \sigma_d: F\to \R,$$
which satisfies the following extra condition: for $1\leq i\leq
d$ there exists $\R$-isomorphism
$$\rho_1: A^{\sigma_1}\otimes\R\simeq M(2,\R),\quad \rho_i:
A^{\sigma_i}\otimes\R\simeq \mathbb H,\quad 2\leq i\leq d,$$
where $\mathbb H$ is the quaternion  algebra over $\R.$ An order
$\sO\subset A$ over $F$ is a subring of $A$ containing 1 which is
a finitely generated $\sO_F-$module generating the algebra $A$
over $F.$ The group of units in $\sO$ of reduced norm 1 is
defined as
$$
\sO^1=\{x\in \sO; \ \rm{Nrd}(x)=1\}.
$$
By Shimura
$\rho_1(\sO^1)\subset Sl_2(\R)$ is a discrete subgroup and
for a torsion free subgroup $\Gamma \subset \sO^1$ of finite index
$\sH/\rho_1\Gamma)$ is a quasi-projective curve, called Shimura
curve. Furthermore, if $A$ is a division algebra
$\sH/\rho_1(\Gamma)$ is projective (see \cite{Shi}, Chapter 9).\\

\begin{remark}\label{ramified}
We will say that over some field extension $F'$ of $F$ the
quaternion algebra splits, if $A_{F'}=A\otimes_FF'\simeq
M(2,F')$. If $F'=F_v$ is the completion of $F$ with respect to a
place $v$ of $F$, one says that $F$ is ramified at $v$, if
$A_v=A_{F_v}$ does not split. As well known, there exists some
$a\in F$ for which $A_{F(\sqrt{a})}$ splits. As explained in
\cite{Vig}, for example, we can choose such $a\in F$ in the
following way:\\[.1cm]
Fix one non-archimedian prime $p^0$ of $\Q$, such that $A$ is
unramified over all places of $F$ lying over $p^0$. Then choose
$a$ such that for all places $v$ of $F$ not lying over $p^0$ the
quaternion algebra $A$ ramifies at $v$ if and only if
$F_v(\sqrt{a})\neq F_v$. Moreover one may assume, that the product
over all conjugates of $a$ is not a square in $\Q$.
\end{remark}
\begin{definition}\label{fuchsian} If $\tilde\Gamma\in PSl_2(\R)$ is a subgroup of
finite index of some $P\rho_1(\sO^1)$, then we call $\tilde\Gamma$
a Fuchsian group derived from a quaternion algebra $A.$
\end{definition}
\begin{theorem}[Takeuchi \cite{Tak}]\label{takeuchi} \quad  Let $\tilde\Gamma\subset
PSl_2(\R)$ be a discrete subgroup such that $\sH/\tilde\Gamma$ is
quasi-projective. Then $\tilde\Gamma$ is derived from a quaternion
algebra $A$ over a totally real number field $F$ with $d$ distinct
embeddings
$$\sigma_1=id,\,\sigma_2,\ldots,\sigma_d: F\>>> \R,$$
with
$$\rho_1: A^{\sigma_1}\otimes\R\simeq M(2,\R),\quad \rho_i:
A^{\sigma_i}\otimes\R\simeq \mathbb H,\quad 2\leq i\leq d$$
if and only if $\Gamma$ satisfies the following conditions:
\begin{enumerate}
\item[(I)] Let $k$ be the field
generated by the set $\rm{tr}(\L)$ over $\Q.$ Then $k$ is an
algebraic number field of finite degree, and $\rm{tr}(\L)$ is
contained in the ring of integers of $k$, $\sO_{k}.$
\item[(II)] Let $\sigma$ be any embedding of $k$ into $\C$ such
that $\sigma\neq{\rm id}_k$. Then $\sigma(\rm{tr}(\L))$ is bounded
in $\C.$
\end{enumerate}
\end{theorem}
In the proof of Theorem \ref{takeuchi} one gets, in fact, $k=F.$
If $A$ is a division algebra, for example if $d>1$, then
$Y=\sH/\tilde\Gamma$ is projective, and it is determined by $A$,
and by the choice of the order $\sO\subset A$ up to finite \'etale coverings.

\begin{assumption}\label{assumption2}
Let $\X_\Q$ be an irreducible $\Q$ variation of Hodge structures
of weight one and width one, and with a maximal Higgs field.
Assume moreover, that $\X_\Q$ is polarized. There are isomorphisms
$$
\psi: \X=\X_\Q\otimes_\Q\C \> \simeq >> \V \oplus \U_1 \mbox{ \ \
and \ \ } \phi:\V \> \simeq >> \L\otimes \T
$$
where $\U_1$ and $\T$ are both unitary, and where $\L$ is a rank
two variation of Hodge structures of weight one and width one,
with a maximal Higgs field. Moreover $\V$, $\U_1$, $\L$, $\T$ and
$\psi$ are defined over some real number field $K$, and $\phi$
over some number field $K'$. We fix an embedding of $K'$ into $\C$
and denote by $k\subset K \subset K' \subset \C$ the field spanned
by $\rm{tr}(\L)$ over $\Q$.
\end{assumption}
\begin{proposition} \label{quatsplit}
Keeping the notations and assumptions made in \ref{assumption2},
replacing $Y$ by an \'etale covering, one may assume that
\begin{enumerate}
\item[i.] $\Gamma_{\L}$ is  derived from a quaternion algebra $A$
over a totally real number field $F$ with $d$ distinct embeddings
$$\sigma_1=id,\,\sigma_2,\ldots,\sigma_d: F\>>> \R.$$
\item[ii.] for $1\leq i\leq d$ there exists $\R$-isomorphism
$$\rho_1: A^{\sigma_1}\otimes\R\simeq M(2,\R),\mbox{ \ and \ }
\rho_i: A^{\sigma_i}\otimes\R\simeq \mathbb H,\mbox{ \ for \ }
2\leq i\leq d.$$
\item[iii.] the representation $\gamma_\L:\pi_1(Y,*)\to \Sl(2,\R)$
defining the local system $\L$ factors like
$$
\pi_1(Y,*)\>\simeq >> \Gamma \subset \rho_1(\sO^1) \>>>
\Sl(2,\R\cap\bar\Q)\subset \Sl(2,\R),
$$
and $Y\simeq \sH/\Gamma$.
\item[iv.] for $a$ as in \ref{ramified}
$F(\sqrt{a})$ is a field of definition for $\L$.
\item[v.] if $\tau_i,\,1\leq i\leq d$ are extension of $\sigma_i$
to $F(\sqrt{a})$, and if $\L_i$ denotes the local system defined
by
$$
\pi_1(Y,*)\>>> \Sl(2,F(\sqrt{a}))\> \tau_i>> \Sl(2,\bar\Q),
$$
then $\L_i$ is a unitary local system, for $i>1$, and $\L_1\simeq \L$.
\item[vi.] up to isomorphism, $\L_i$ does not depend on
the extension $\tau_i$ chosen.
\end{enumerate}
\end{proposition}

\begin{proof} \ i) and ii): By
Corollary \ref{split6}, iii), $\Gamma_{\L}$ satisfies Condition
(I) in Theorem \ref{takeuchi}. So, we only have to verify Condition
(II) for $\Gamma_\L$. Let $\sigma$ be an embedding of $k$ into
$\C$ which is not the identity, and let $\tilde \sigma: K'\to \C$
be any extension of $\sigma$.

By \ref{simpson3} $\psi^{-1}\V^{\tilde\sigma}$ is a sub variation
of Hodge structures of $\X$, hence of width zero or one. On the
other hand, $\V^{\tilde\sigma}$ is isomorphic to
$\L^{\tilde\sigma}\otimes\T^{\tilde\sigma}$. Both factors are
variations of Hodge structures, hence at least one of them has a
trivial Higgs field.

Assume both have a trivial Higgs field, hence $\V^{\tilde\sigma}$
as well. By \ref{width}, a), the composite
$$
\psi^{-1}\V^{\tilde\sigma} \>>> \X \> \psi >> \V\oplus\U_1 \>>> \V
$$
has to be zero. Hence $\V^{\tilde\sigma}$ is a sublocal system of
the unitary system $\U_1$, hence unitary itself. The $\bar\Q$
isomorphism $\phi: \V \> \simeq >> \L \otimes \T$ induces an
isomorphism
$$
\phi^\otimes : \bigotimes^{g_0} \V^{\tilde\sigma} \>>> \big(
\bigotimes^{g_0} \L^{\tilde\sigma}\big) \otimes \big(
\bigotimes^{g_0} \T^{\tilde\sigma}\big).
$$
The right hand side contains $S^{g_0} (\L^{\tilde\sigma})$ as a
direct factor, hence $S^{g_0} (\L^{\tilde\sigma})$ is unitary, as
well as $\L^{\tilde\sigma}$. So
$\tr(\L^{\tilde\sigma})=\sigma(\tr(\L))$ is bounded in this case.

If the Higgs field of $\V^{\tilde\sigma}$ is non trivial, it is
generically maximal. This implies that the composite
$$
\V^{\tilde\sigma} \> \simeq >> \X \> >> \U_1
$$
is zero. Hence $\V^{\tilde\sigma}\simeq \V$. If the Higgs field of
$\L^{\tilde\sigma}$ is an isomorphism, by \ref{semistableii}
replacing $Y$ by an \'etale covering, $\L^{\tilde \sigma}\simeq
\L$. Hence up to conjugation the representations
$\gamma_{\L^{\tilde\sigma}}$ and $\gamma_{\L}$ coincide and for
all $\eta\in\pi_1(Y,*)$
$$
\tr(\gamma_{\L^{\tilde\sigma}}(\eta)) = \tr(\gamma_\L(\eta)).
$$
So $\sigma$ is the identity, a contradiction.\\

It remains to consider the case that $\V^{\tilde\sigma}\simeq \V$ and
that $\L^{\tilde\sigma}$ is concentrated in degree $0,0$.

For $g_0$ even,  one has a $\bar\Q-$isomorphism
$$ \wedge^{g_0}\phi: \wedge^{g_0}\V\simeq S^{g_0}(\L)\oplus {\mathbb
S}_{\{2,\ldots,2\}}(\T) \oplus \bigoplus_{c=1}^{\frac{g_0}{2}-1}
S^{2c}(\L)\otimes {\mathbb S}_{\lambda_{2c}}(\T),
$$
where ${\mathbb S}_{\{2,\ldots,2\}}(\T)$ is of width zero, where
$S^{g_0}\L$ has a maximal Higgs field of width $g_0$, and where
all other factors have a maximal Higgs field of width between $2$
and $g_0-2$. Let $K$ denote the field of definition $\Gamma_\L.$
Then $K\supset k$ is a finite extension of $k.$ Let $\sigma$ be an
embedding of $k$ into $\C$ which is not identity, and let $\tilde
\sigma: K\to \C$ be an extension of $\sigma.$ Via the
isomorphisms $\wedge^{g_0}\phi$ and
$\wedge^{g_0}\phi^{\tilde\sigma}$ we obtain an embedding
$$
S^{g_0}\L^{\tilde \sigma}\>>> S^{g_0}(\L)\oplus {\mathbb
S}_{\{2,\ldots,2\}}(\T) \oplus \bigoplus_{c=1}^{\frac{g_0}{2}-1}
S^{2c}(\L)\otimes {\mathbb S}_{\lambda_{2c}}(\T).
$$
The projection of $S^{g_0}\L^{\tilde \sigma}$ into $S^{g_0}\L$
must be zero, for otherwise, we would get an isomorphism
$S^{g_0}\L^{\tilde \sigma}\simeq S^{g_0}\L.$ By Corollary
\ref{simpson3} $\L^{\tilde \sigma}$ is a sub variation of Hodge
structures, hence it has a maximal Higgs field.

The projection
$$S^{g_0}\L^{\tilde \sigma}\>>>
\bigoplus_{c=1}^{\frac{g_0}{2}-1} S^{2c}(\L)\otimes {\mathbb
S}_{\lambda_{2c}}(\T)
$$ must be also zero, for otherwise, by
applying again Corollary \ref{simpson3} to $S^{g_0}\L^{\tilde
\sigma}$ one would find $\L^{\tilde \sigma}$ to have a maximal
Higgs field, hence $S^{g_0}\L^{\tilde \sigma}$ to have a
maximal Higgs field of width $g_0.$ But, then it can not be
embedded in a local system of width $<g_0.$

Thus, the projection
$$ S^{g_0}\L^{\tilde \sigma}\>>> S_{\{2,2,\ldots,2\}}\T$$
is an embedding. This implies that $\L^{\tilde \sigma}$ is
unitary. In particular, again
$\rm{tr}(\L^{\tilde\sigma})=\sigma(\rm{tr}(\L))$ is bounded in
$\C.$

Finally, the assumption $\V^{\tilde\sigma}\simeq \V$ and
$\L^{\tilde\sigma}$ unitary does not allow $g_0={\rm rank}(\V)$
to be odd:

The $\bar\Q$ isomorphism $\phi: \V \simeq \L\otimes\T,$
induces a $\bar\Q$ isomorphism
$$ \wedge^{g_0}\phi: \wedge^{g_0}\V\simeq \bigoplus_{c=0}^{\frac{g_0-1}{2}}
S^{2c+1}(\L)\otimes {\mathbb S}_{\lambda_{2c}}(\T)
$$
(see \ref{tensor}). The left hand side contains a local subsystem
isomorphic to $S^{g_0}(\L^{\tilde\sigma})$, hence with a trivial Higgs field,
whereas the right hand side only contains factors of width $>0$,
with a maximal Higgs field, a contradiction.\\
\ \\
Applying \ref{takeuchi} we obtain a quaternion algebra $A$
satisfying i), ii) and the first part of iii). By \ref{quotient}
one has $Y\simeq \sH/\Gamma$.

For iv) we recall that by the choice of $a$ the quaternion
algebra $A$ splits over $F(\sqrt{a})$. So v) follow from i) and
ii).

To see that $\L_i$ is independent of the extension of $\sigma_i$
to $\tau_i: F(\sqrt{a}) \to \bar\Q$ it is sufficient to show vi)
for $i=1$. Let $\bar\L$ denote the local system obtained by
composing the representation with the involution on
$F(\sqrt{a})$. Then both, $\L$ and $\bar\L$ have a maximal Higgs
field, hence by \ref{semistableii}, c), their Higgs fields differ
at most by the product with a two torsion element in ${\rm
Pic}^0(Y)$. Replacing $Y$ by an \'etale covering, we may assume
both to be isomorphic.
\end{proof}

Given a quaternion algebra $A$ as in \ref{quatsplit}, i) and ii)
allows to construct certain families of abelian varieties. To this
aim we need some well known properties of quaternion algebras $A$
defined over number fields $F$. Let us fix a subfield $L$ of $F$.
\begin{notations}\label{notations}
Let us write $\delta=[L:\Q]$, $\delta'=[F:L]$ and
$$
\beta_1={\rm id}_L, \beta_2, \ldots , \beta_\delta : L \>>> \C
$$
for the different embeddings. We renumber the embeddings
$\sigma_i:F \to \C$ in such a way, that
$$
\sigma_i|_L = \beta_\nu \mbox{ \ \ for \ \ } (\nu-1)\delta' < i
\leq \nu \delta'.
$$
\end{notations}

Recall, that the corestriction $\Cor_{F/L}(A)$ is defined (see
\cite{Vig}, p. 10) as the subalgebra of ${\rm Gal}(\bar\Q/L)$
invariant elements of
$$
\bigotimes_{i=1}^{\delta'} A^{\sigma_i} =
\bigotimes_{i=1}^{\delta'} A\otimes_{F,\sigma_i}\bar\Q.
$$
\begin{lemma}\label{corestriction}
Let $A$ be a quaternion division algebra defined over a totally real number
field $F$, of degree $d$ over $\Q$. Assume that $A$
is ramified at all infinite places of $F$ except one. For some
subfield $L$ of $F$ let $D_L=\Cor_{F/L}(A)$ be the corestriction of $A$
to $L$. Finally let $a\in F$ be an element, as defined in \ref{ramified},
and
$$
b=a\cdot\sigma_2(a)\cdot \cdots \cdot \sigma_{\delta'}(a) \in L.
$$
\begin{enumerate}
\item[a.] If $L=\Q$, i.e. if $d=\delta'$, then either
\begin{enumerate}
\item[i.] $D_\Q\simeq M(2^{d},\Q)$, and $d$ is odd, or
\item[ii.] $D_\Q\not\simeq M(2^{d},\Q)$. Then
$$
D_\Q\simeq M(2^{d},\Q(\sqrt{b})).
$$
$\Q(\sqrt{b})$ is a quadratic extension of $\Q$, real if and only if $d$
is odd.
\end{enumerate}
\item[b.] If $L\neq \Q$, then $D_L\not\simeq M(2^{\delta'},L)$, and
\begin{enumerate}
\item[i.] $L(\sqrt{b})$ is an imaginary quadratic extension of
$L$.
\item[ii.] $ D_L\otimes_LL(\sqrt{b})\simeq M(2^{\delta'},L(\sqrt{b}))$.
\end{enumerate}
\end{enumerate}
In a), ii), or in b),
choosing an embedding $L(\sqrt{b}) \to M(2,L)$, one obtains an
embedding
$$
D_L \>>> M(2^{d+1},L).
$$
\end{lemma}
\begin{proof}
For $\delta=[L:\Q] \geq 1$, choose $\delta$ different embeddings
$\beta_\nu: L \to \bar\Q$, corresponding to infinite places $v_1,
\ldots , v_\delta$. We may assume that $\beta_1$ extends to the
embedding $\sigma_1$ of $F$. Hence $A$ is ramified over
$\delta'-1$ extensions of $v_1$ to $F$, and over all $\delta'$
extensions of $v_\nu$ to $F$, for $\nu \neq 2$. Writing $L_v$ for
the completion of $L$ at $v$, one has
$$
D_{v_\nu}=\Cor_{F/L}A\otimes_LL_{v_\nu}=\left\{
\begin{array}{ll}
M(2,\R)\otimes \bigotimes^{\delta'-1}\HH & \mbox{ for } \nu = 1\\
\bigotimes^{\delta'}\HH & \mbox{ for } \nu \neq 1
\end{array} \right. .
$$
Recall that the $r$-fold tensor product of $\HH$ is isomorphic to
$M(2^r,\R)$ if and only if $r$ is even. By our choice of $a$ and
$b$ this holds true, if and only if $L_v(\sqrt{b})=L_v$. In fact,
the image of $b$ in $L_{v_\nu}$ has the sign $(-1)^{\delta'-1}$,
for $\nu=1$ and $(-1)^{\delta'}$ otherwise.

In particular $D_L\simeq M(2^{\delta'},L)$ can only hold true for
$L=\Q$ and $d=\delta'$ odd.

For $L=\Q$, one also finds $b>0$ if and only if $d$ is odd.

For all but finitely many non-archimedian places $v$ of $L$, in
particular for those dominating the prime $p^0$ in \ref{ramified},
and for the completion $L_v$ with respect to $v$, one has
$$
D_v=\Cor_{F/L}A\otimes_LL_{v}=M(2^{\delta'},L_v).
$$
If this is not the case, consider the extension
$L'_v=L(\sqrt{b})\otimes_{L}L_v$ of $L_v$. One finds
$$
D_v\otimes_{L_v}L'_v = M(2^{\delta'},L'_v).
$$
In fact, let $v_1,\ldots,v_\ell$ be the places of $F$, lying over
$v$, and let $F_1,\ldots,F_\ell$ be the corresponding local
fields. Then
$$
D_v = \bigotimes_{i=1}^\ell \Cor_{F_i/L_v} (A\otimes_FF_i),
$$
and it is sufficient to show that $D_i=\Cor_{F_i/L_v}
(A\otimes_FF_i)$ splits over $L'_v$. If $L'_v$ is a subfield of
$F_i$
$$
D_i\otimes L'_v= \Cor_{F_i/L_v} (A\otimes_FF_i)^{\otimes 2}
$$
splits, since $(A\otimes_FF_i)^{\otimes 2}$ does. The same holds
true, if $L'_v$ is not a field. If $L'_v$ is a field, not
contained in $F_i$, then $F'_i=F_i\otimes_{L_v}L'_v$ is a field
extension of $F_i$ of degree two, and
$$
D_i\otimes L'_v= \Cor_{F'_i/L_v} (A\otimes_FF'_i)
$$
splits again, since $(A\otimes_FF'_i)$ does.

By \cite{Wei}, Chapter XI, \S 2, Theorem 2 (p. 206),
$$
D_{L(\sqrt{b})}=D\otimes_LL(\sqrt{b})=M(2^{\delta'},L(\sqrt{b})).
$$
\end{proof}

Choose again an order $\sO$ in $A$, and let $\sO^1$ be the group
of units in $\sO$ of reduced norm $1$. For any discrete torsion free
subgroup $\tilde\Gamma \subset P\rho_1(\sO^1)$ with preimage $\Gamma$
in $\sO^1\subset\Sl_2(\R)$ the diagonal embedding
\begin{equation*}
\Gamma \>>> \sO^1 \>>> \bigotimes_{i=1}^{\delta'} A^{\sigma_i}
\end{equation*}
induces an embedding
\begin{equation}\label{repr}
\Gamma \>>> \sO^1 \>>> D_L=\Cor_{F/L}A.
\end{equation}
\begin{construction}\label{constrq} \
For $L=\Q$ the morphism (\ref{repr}) and \ref{corestriction}, a),
give a morphism
$$
\Gamma \subset D=\Cor_{F/\Q}A \subset
D\otimes_\Q\Q(\sqrt{b})=M(2^d,\Q(\sqrt{b})) \subset
M(2^{d+\epsilon},\Q)
$$
for $\epsilon = 0$ or $1$, where $b\in \Q$ is either a square, or
as defined in \ref{quatsplit}. One obtains a representation
$$
\eta:\Gamma \>>> \Gl(2^d,\Q(\sqrt{b})) \>>>
\Gl(2^{d+\epsilon},\Q).
$$
If $\epsilon=0$, the degree $d$ must be odd. Over $\R$ one has
\begin{equation}\label{structure}
D\otimes_\Q\R \simeq M(2,\R)\otimes \HH \otimes \cdots \otimes
\HH.
\end{equation}
The $\Q$ algebraic group $G:=\{ x\in D^*; \ {\rm Nrd}(x)=x{\bar
x}=1\}$ is $\Q$ simple and by (\ref{structure}) it is a $\Q$-form
of the $\R$ algebraic group
$$
G(\R)\simeq \Sl(2,\R)\times {\rm SU}(2)\times \cdots \times {\rm
SU}(2).
$$
Projection to the first factor, gives a representation of $\Gamma$
in $\Sl(2,\R)$, hence a quotient $Y=\sH/\Gamma$ with
$\Gamma=\pi_1(Y,*)$.

Let us denote by $\V_\Q$ or by $\X_\Q$ the $\Q$ local system on $Y$
induced by $\eta$. If we want to underline,
that the local systems are determined by $A$ we also write
$\V_{A\Q}$ and $\X_{A\Q}$, respectively.
\end{construction}
\begin{lemma}\label{constrq1}
Keeping the assumptions and notations from \ref{constrq} one
finds:
\begin{enumerate}
\item[a.] $$
\dim(\End(X_{A,\Q}))= \left\{ \begin{array}{ll}
1 & \mbox{for } \epsilon=0\\
4 & \mbox{for } \epsilon=1 \end{array}\right. .
$$
\item[b.] For $\epsilon=1$ one has
$$
\dim(H^0(Y,\bigwedge^2(X_{A,\Q})))= \left\{ \begin{array}{ll}
3 & \mbox{for } d \mbox{ odd}\\
1 & \mbox{for } d \mbox{ even} \end{array}\right. .
$$
\end{enumerate}
\end{lemma}
\begin{proof}
Consider for $\epsilon'=2^\epsilon$
$$
\X=\X_{A\Q}\otimes \C= \L_1\otimes \cdots \otimes \L_d \otimes
\C^{\epsilon'}
$$
where for $\tilde\sigma\in {\rm Gal}(\bar\Q/\Q)$ the local system
$\L_i^{\tilde\sigma}$ has a maximal Higgs field if and only if
$\tilde\sigma|_F=\sigma_i^{-1}$. Otherwise this local system is
unitary and of pure bidegree $0,0$.

The determinant of each $\L_i$ is $\C$, hence $\End(\X)$ contains
$\C^{\epsilon'}\otimes \C^{\epsilon'}$ as a direct factor. Then
\begin{equation}\label{sections3}
\dim_\Q(\End(\X_{A\Q}))=\dim_\C(|End(\X)) \geq 4^\epsilon.
\end{equation}
One has
$$
\End(\X_\Q)=H^0(Y,\END(\X_\Q)) \simeq
H^0(Y,\bigwedge^2(\X_\Q))\oplus H^0(Y,S^2(\X_\Q)).
$$
By \ref{tensor1}
\begin{gather*}
H^0(Y,\bigwedge^2(\X))= H^0(Y,S^2(\L_{2}\otimes \cdots \otimes
\L_{d}
\otimes \C^{\epsilon'})) \mbox{ \ \ \ and}\\
H^0(Y,S^2(\X))= H^0(Y,\bigwedge^2(\L_{2} \otimes \cdots \otimes
\L_{d} \otimes \C^{\epsilon'})).
\end{gather*}
Since $\End(\X_\Q)$ is invariant under ${\rm Gal}(\bar\Q/\Q)$, for
$\tilde\sigma$ with $\tilde\sigma|_F=\sigma_2$ it is for $d>1$
contained in the direct sum of
\begin{gather*}
H^0(Y,S^2(\L^{\tilde\sigma}_{2}\otimes \cdots \otimes
\L^{\tilde\sigma}_{d} \otimes
\C^{\epsilon'}))=H^0(Y,\bigwedge^2(\L^{\tilde\sigma}_{3}\otimes
\cdots \otimes \L^{\tilde\sigma}_{d}
\otimes \C^{\epsilon'})) \mbox{ \ \ \ and}\\
H^0(Y,\bigwedge^2(\L^{\tilde\sigma}_{2}\otimes \cdots \otimes
\L^{\tilde\sigma}_{d} \otimes
\C^{\epsilon'}))=H^0(Y,S^2(\L^{\tilde\sigma}_{3}\otimes \cdots
\otimes \L^{\tilde\sigma}_{d} \otimes \C^{\epsilon'})).
\end{gather*}
Repeating this game we find
\begin{gather}
H^0(Y,\bigwedge^2(\X_\Q))\subset \left\{
\begin{array}{ll}
H^0(Y,S^2(\C^{\epsilon'}))& \mbox{for } d \mbox{ odd}\\
H^0(Y,\bigwedge^2(\C^{\epsilon'}))& \mbox{for } d \mbox{ even}
\end{array}\right. \label{sections1} \\
H^0(Y,S^2(\X_\Q))\subset \left\{
\begin{array}{ll}
H^0(Y,\bigwedge^2(\C^{\epsilon'}))& \mbox{for } d \mbox{ odd}\\
H^0(Y,S^2(\C^{\epsilon'}))& \mbox{for } d \mbox{ even}
\label{sections2}
\end{array}\right. .
\end{gather}
For $\epsilon'=1$ we obtain that $\End(\X_\Q)$ is a most one dimensional
and for $\epsilon'=2$ we find
$$
\dim_\Q(H^0(Y,\bigwedge^2(\X_\Q))) \leq 3 \mbox{ \ \ and \ \ }
\dim_\Q(H^0(Y,S^2(\X_\Q))) \leq 1
$$
or vice versa. Comparing this with (\ref{sections3})
one obtains \ref{constrq1} i) and ii).
\end{proof}
\begin{lemma}\label{constrq2}
Given a quaternion division algebra $A$, as in \ref{quatsplit} i)
and ii), there exists a smooth family of abelian varieties $f:X_A\to
Y$ with $R^1f_*\Q_{X_A}=\X_{A\Q}$. Moreover, the special Mumford-Tate
group $\Hg$ of the general fibre of $f$ is the same as the group
$G$ in \ref{constrq}.
\end{lemma}
\begin{proof} (see \cite{Mum2})
The group $G$ in \ref{constrq} and the representation
$$
G\>>> D^*\to \Gl(2^{d+\epsilon},\Q)
$$
are $\Q$ forms of an $\R$-representation
$$
\Sl(2,\R)\times {\rm SU}(2)^{\times (d-1)}\>>> \Sl(2,\R)\times
{\rm SO}(2^{d-1})\>>> \Gl(2^{d+\epsilon},\R).
$$
The group in the middle acts on $\R^2\times \R^{2^{d-1}}$. Over
$\R,$ this representation leaves a unique non degenerate
symplectic form $< \ , \ >$ on $\R^{2^{d}}$ invariant, the tensor
product of the $\Sl(2,\R)$ invariant symplectic form on $\R^2$
with the ${\rm SO}(2^{d-1})$ invariant Hermitian form.

Hence for $\epsilon=0$ and $V=\Q^{2^{d}}$ there is a unique
symplectic form $Q$ on $V$, invariant under $\Gamma\subset G$.

For $\epsilon=1$, one chooses $V=\Q(\sqrt{b})^{2^d}$. Again one has
a unique $\Q(\sqrt{b})$ valued symplectic form on $V$. Regarding
$V$ as a $\Q$ vector space, the trace $\Q(\sqrt{b}) \to \Q$ gives
a $\Q$ valued symplectic form $Q$, again invariant under $\Gamma
\subset G$.

Note that $\Gamma$ is the group of units of an order $\mathcal O$
in $A$. Hence $\Gamma$ leaves a $\Z$-module $L\subset V$ of rank
$\dim V$ invariant. For some submodule $H\subset L$ of the form
$H=mL$, for $m\gg 0$, one has $Q(H\times H)\subset \Z$. Obviously
$\Gamma$ leaves $H$ again invariant. So one obtains a
representation
$$
\Gamma\>>> {\rm Sp}(H,Q)\otimes\Q.
$$
Finally let
$$\phi_0:T=\{ z \in \C; \ |z|=1 \} \>>> \Sl(2,\R)\times {\rm SO}(2^{d-1})
\subset {\rm Sp}(H, Q)\otimes \R
$$
be the homomorphism defined by
$$
e^{i\theta}\mapsto \left[ \begin{array}{cc}
 \cos\theta & \sin\theta\\
-\sin\theta & \cos\theta
\end{array} \right]\times {\rm I}_{2^{d-1}}.
$$
$J_0=\phi_0(i)$ defines a complex structure on $H\otimes\R,$ and
$$
Q(x, J_0x)>0, \mbox{ \ \ for all \ \ } x\in H.
$$
The image of $G$ in ${\rm Sp}(H,Q)\otimes \R$ is normalized by
$\phi_0(T),$ i.e. for all $g\in G$ one has
$$
g\phi_0(T)g^{-1}=\phi_0(T).
$$
So $\X_{A\Q}$ defines a smooth family
of abelian varieties $f: X_A\to Y=\sH/\Gamma.$

By the construction this family reaches the Arakelov bound and
$\X_{A\Q}$ has no unitary part. By Lemma \ref{hg}, c), one knows that
$$
G_0^{\rm Mon} = {\rm Hg}(R^1f_*\Q_{X_A}).
$$
On the other hand, $G_0^{\rm Mon}$ is contained in the image of
$G$ in ${\rm Sp}(H,\Q)\otimes\Q$. Since
$$
\X_{A\C}=\L_{1\C}\otimes\L_{2\C}\otimes\cdots\otimes\L_{d\C}\otimes
\C^{2^\epsilon}
$$
and since all factors are Zariski dense in $\Sl(2,\C)$ one finds
that
$$
G_{0\C}^{\rm Mon}=\Sl(2,\C)^{\times d}=G_\C,
$$
hence
$$
G_0^{\rm Mon} = {\rm Hg}(R^1f_*\Q_{X_A}) = G.
$$
\end{proof}
Let us remark, that in the proof of Theorem \ref{shimura} in
Section \ref{sectcorestriction} we will see, that the families
$f:X_A \to Y$ in \ref{constrq2} are unique up to isogenies, and up
to replacing $Y$ by \'etale coverings, and that they belong to one
of the examples described in \ref{mumfordshimura}.

\begin{construction}\label{constrl}
If $L\neq \Q$ choose $b$ as in \ref{corestriction}. The morphism
(\ref{repr}) and \ref{corestriction}, b), give a map
$$
\Gamma \subset D_L=\Cor_{F/L}A \subset
D_L\otimes_LL(\sqrt{b})=M(2^{\delta'},L(\sqrt{b})) \subset
M(2^{\delta'+1},L),
$$
inducing a representation $\Gamma \to \Gl(2^{\delta'+1},L)$,
hence an $L$ local system $\V_L$ on $Y=\sH/\Gamma$.

An embedding $L\subset M(\delta,\Q)$ gives rise to
$$
\Gamma \subset D_L=\Cor_{F/L}A \subset M(2^{\delta'+1},L) \subset
M(\delta 2^{\delta'+1},\Q),
$$
hence to a $\Q$ local system $\X_\Q=\X_{A,L;\Q}$.

In different terms, choose extensions $\tilde\beta_\nu$ of
$\beta_\nu$ to $\bar\Q$. For  $\V_{\bar\Q}=\V_L\otimes_L\bar\Q,$
the $\bar\Q$ local system
$$
X_{\bar\Q}=\X_{A,L;\bar\Q}=
\V_{\bar\Q} \oplus \V_{\bar\Q}^{\tilde\beta_2} \oplus
\cdots \oplus \V_{\bar\Q}^{\tilde\beta_\delta}
$$
is invariant under ${\rm Gal}(\bar\Q/\Q)$,
hence defined over $\Q$.
\end{construction}

\begin{remark}
Consider any family $X\to Y$ of abelian varieties, with a geometrically simple
generic fibre. If $\X_{A,L;\Q}$ is an irreducible component of
$R^1f_*\Q_X$, all irreducible components of $R^1f_*\Q_X$ are
isomorphic to $\X_{A,L;\Q}$. As in \cite{Del}, p. 55, for
$\Delta=\End(\X_{A,L;\Q})$ one finds
$$
R^1f_*\Q_X \simeq \X_{A,L;\Q} \otimes_\Delta {\rm
Hom}(\X_{A,L;\Q},R^1f_*\Q_X),
$$
and for some $m$
$$
\End(R^1f_*\Q_X) \simeq M(m,\Delta).
$$
In \cite{Sai}, Section 9, one finds examples showing that all
$m>0$ occur.
\end{remark}

\section{The proof of Theorems \ref{shimura} and \ref{shimura2}}
\label{sectcorestriction}

In order to prove Theorems \ref{shimura} and \ref{shimura2} we
will show, that the local subsystem $\X_\Q$ in \ref{assumption2}
is for some $L\subset F$ isomorphic to the one constructed in
\ref{constrq} or \ref{constrl}.

Let us consider the subgroup $H$ of ${\rm Gal}(\bar\Q/\Q)$ of all
$\beta$ with $(\psi^{-1}\V)^\beta=\psi^{-1}\V$, and let $L$ denote
the field of invariants under $H$. So $\V=\V_L\otimes_L\C$.

\begin{proposition}\label{quatsplit2} Let us keep the assumptions
made in \ref{assumption2} and use the notations introduced in
\ref{quatsplit}. Replacing $Y$ by a finite \'etale covering, the field
of invariants $L$ under $H$ is a subfield of $F$. Using the notations
introduced in \ref{notations} for such a subfield,
there exists a decomposition $\V_L
\simeq \L_{1L} \otimes \cdots \otimes \L_{\delta' L} \otimes
\T'_L$ with:
\begin{enumerate}
\item[i.] For $\beta\in {\rm Gal}(\bar\Q/L)$ and $i\leq \delta'$ one has
$\L^\beta \simeq \L_i$, if and only if $\beta|_F=\sigma_i$.
\item[ii.] For $\beta \in {\rm Gal}(\bar\Q/L)$ the Higgs field of
$(\L_1\otimes \cdots \otimes \L_{\delta'})^\beta$ is maximal.
\item[iii.] For $\beta \in {\rm Gal}(\bar\Q/\Q)$ with $\beta|_L\neq {\rm
id}_L$ the local system $(\L_1\otimes \cdots \otimes
\L_{\delta'})^\beta$ is unitary.
\item[iv.] For $\beta\in {\rm Gal}(\bar\Q/L)$ the local system
${\T'}^\beta$ is a unitary.
\end{enumerate}
\end{proposition}
\begin{proof} Replacing $Y$ by an \'etale covering, we are allowed to
apply \ref{quatsplit}. In particular we have the
rank $2$ local systems $\L_1, \ldots , \L_d$, defined there.
Consider any decomposition
$\V\simeq \L_1 \otimes \cdots \otimes \L_r \otimes \T_r$ with:
\begin{enumerate}
\item[i'.] If for $\beta\in {\rm Gal}(\bar\Q/\Q)$ one has
$\beta|_F=\sigma_i$, with $i\in \{1, \ldots , r\}$, then
$\L^\beta \simeq \L_i$.
\end{enumerate}
For $r=1$, \ref{semistableii} gives a decomposition
$\V=\L\otimes \T$. Write again $\L_1=\L$ and $\T_1=\T$.
By \ref{quatsplit}, iv), the local system is defined over
$F(\sqrt{a})$ and by \ref{quatsplit} vi),
$\L_1^\beta\simeq \L_1$ if the restriction of $\beta$ to $F$ is
$\sigma_1={\rm id}_F$. Hence i') holds true for this decomposition.

Consider for $r\geq 1$ a decomposition satisfying i'). \\
\ \\
{\bf Step 1.} If for some $\beta'\in {\rm Gal}(\bar\Q/\Q)$
and for $i\in \{1, \ldots , r\}$ one has $\L^{\beta'} \simeq \L_i$, then
necessarily $\beta'|_F=\sigma_i$.\\
\ \\
In fact, let $\beta \in {\rm Gal}(\bar\Q)$ be an automorphism with
$\beta_F=\sigma_i$. Then $\L^{\beta^{-1}\circ\beta'}\simeq \L$,
and \ref{quatsplit}, v), implies that $\beta^{-1}\circ\beta'|_F={\rm id}_F$.\\
\ \\
{\bf Step 2.} There exists no $\tau\in {\rm Gal}(\bar\Q/\Q)$ with
$\L_1^\tau\otimes\cdots\otimes\L_r^\tau$ not unitary and with
$\tau|_F\neq \sigma_i$ for $i=1, \ldots ,r$.\\
\ \\
Assume the contrary. Renumbering the embeddings of $F \to \R$ one
may assume that $\tau|_F=\sigma_{r+1}$. Recall that
by \ref{simpson3} $\L_i^\tau$ is a variation of Hodge structures
of rank $2$. It either is of width zero, hence unitary, or of
width one, hence with maximal Higgs field. By assumption there exists some
$i<r+1$ for which $\L_i^\tau$ has a maximal Higgs field. Choose
$\beta\in {\rm Gal}(\bar\Q/\Q)$ with $\beta|_F=\sigma_i$. Then
$\L^{\beta\circ\tau}=\L_i^\tau$ has a maximal Higgs field. \ref{quatsplit},
v), implies $\beta\circ\tau|_F={\rm id}_F$, a contradiction.\\
\ \\
{\bf Step 3.} Assume there exists some $\tau\in {\rm Gal}(\bar\Q/\Q)$ with
$\tau|_F\neq \sigma_i$ for $i=1, \ldots ,r$, with $\V^\tau$ not
unitary, but with $\L_1^\tau \otimes\cdots\otimes\L_r^\tau$ unitary.
Then (renumbering the embeddings $F\to \R$, if necessary)
one finds a decomposition with $r+1$ factors, satisfying
again i').\\
\ \\
$\L_i^\tau$ is unitary for $i=1, \ldots , r$. By
\ref{semistableii} and by \ref{addendum} over some \'etale
covering of $Y$ we find a splitting
$\T_r^\tau\simeq \L\otimes
\T''$, with
$$
\V^\tau\simeq (\L_1 \otimes \cdots \otimes
\L_r)^\beta \otimes \L\otimes \T'').
$$
Apply $\tau^{-1}$. Then one has
$$
\V\simeq \L_1 \otimes \cdots
\otimes \L_r \otimes\L^{\tau^{-1}} \otimes \T_{r+1}.
$$
Since $\L_1$ has maximal Higgs field, $\L_{r+1}:=\L^{\tau^{-1}}$
must be unitary, as well as
$\T_{r+1}$. Applying any extension $\tau_i$ of
$\sigma_i^{-1}$ for $i \leq r$, one finds $\L_{r+1}^{\tau_i}$ to be
unitary, since otherwise there would be two factors with a maximal
Higgs field, $\L_i^{\tau_i}$ and  $\L_{r+1}^{\tau_i}$.

So $\tau|_F$ must be one of the remaining $\sigma_j$, and
renumbering we may assume $\tau|_F=\sigma_{r+1}$.\\
\ \\
{\bf Step 4.} Assume we have found a decomposition
as in i'), and of maximal possible length.
Then for all $\tau\in {\rm Gal}(\bar\Q/\Q)$ with
$\tau|_F\neq \sigma_i$ for $i=1, \ldots ,r$ the local system
$$
\V^\tau \simeq\L_1^\tau
\otimes\cdots\otimes\L_r^\tau\otimes\T^\tau_r
$$
is unitary. For those $\tau$ one has $(\psi^{-1}\V)^\tau\neq
\psi^{-1}\V$. On the other hand, for all $\beta$ with
$\beta|_F=\sigma_i$ with $1\leq i \leq r$ the local system
$\V^\beta$ has a maximal Higgs field, hence
$(\psi^{-1}\V)^\beta=\psi^{-1}\V$. So
$$
H=\{ \beta \in {\rm Gal}(\bar\Q/\Q); \ \beta|_F=\sigma_i \mbox{
with } 1\leq i \leq r \}
$$
and $L$ as the field of invariants under $H$ is contained in $F$.
Using the notations introduced in \ref{notations} for such subfields,
one finds $r=\delta'$ and
$\L_1^\beta\otimes \cdots \otimes \L_{\delta'}^\beta$ has a
maximal Higgs field, for all $\beta \in H$. This in turn implies
that ${\T'}^\beta$ is unitary for those $\beta$.
\end{proof}

\begin{theorem}\label{quatsplit3}
Let us keep the assumption made in \ref{assumption2} and use the
notations introduced in \ref{notations}. Replacing $Y$ by an
\'etale covering, there exists some $\epsilon'>0$ and a
decomposition
\begin{equation}\label{dec}
\psi: \X \> \simeq >> \bigoplus_{\nu=1}^{\delta}\big(
\bigotimes_{i=(\nu-1)\delta'}^{\nu\delta'} \L_i \big)^{\oplus
\epsilon'}
\end{equation}
such that:
\begin{enumerate}
\item[a.] For $\beta\in {\rm Gal}(\bar\Q/\Q)$ the local system
$\L_i^{\beta^{-1}}$ has a maximal Higgs field if and only if
$\beta|_F=\sigma_i$. Moreover $\L^\beta=\L_i$ in this case.
\item[b.] The direct sum in (\ref{dec}) is orthogonal with respect
to the polarization.
\item[c.] If the local subsystems $\psi^{-1}
\L_1\otimes \ldots \otimes \L_{\delta'}$ of $\X$ are defined over
$L$ then $\epsilon'=1$, $L=\Q$ and $[F:\Q]$ is odd.
\item[d.] If $\psi^{-1} \L_1\otimes \ldots \otimes \L_{\delta'}
\subset \X$ is not defined over $L$ choose $b$ to be the element
defined in \ref{quatsplit} and $\iota\in {\rm Gal}(\bar\Q/L)$ with
$\iota(\sqrt{b})=-\sqrt{b}$. Then $\epsilon'=2$, the direct factor
$\psi^{-1} \L_1\otimes \ldots \otimes \L_{\delta'}\otimes \C^2$
in (\ref{dec}) is defined over $L$ and it decomposes over
$L(\sqrt{b})$ like
$$
\psi^{-1} \L_1\otimes \ldots \otimes \L_{\delta'}\oplus
(\psi^{-1} \L_1\otimes \ldots \otimes \L_{\delta'})^\iota \subset
\X.
$$
\item[e.] $\L_1\otimes \cdots \otimes \L_{\delta'}$ is irreducible
as a $\C$ local system.
\end{enumerate}
\end{theorem}
\begin{proof}
Using the notations from \ref{quatsplit2} let us define
$\L_i=\L^{\tilde\sigma_i}$, where $\tilde\sigma_i$ is any
extension of $\sigma_i$ to $\bar\Q$. Obviously, fixing any
extension $\tilde\beta_\nu$ of $\beta_\nu$ one has
$$
\V^{\tilde\beta_\nu} = \L_{(\nu-1)\delta'+1}\otimes\cdots\otimes
\L_{\nu\delta'} \otimes {\T'}^{\tilde\beta_\nu}.
$$
$\V$ has a maximal Higgs field, whereas
$\bigoplus_{\nu=2}^{\delta}\V^{\tilde\beta_\nu}$ is unitary. Hence
their intersection is zero. Applying $\tilde\beta_\nu$ one
obtains the same for the intersection of $\V^{\tilde\beta_\nu}$
and $\bigoplus_{\mu=1,\mu\neq\nu}^{\delta}\V^{\tilde\beta_\mu}$.
So
$$
\psi^{-1}\big(\bigoplus_{\nu=1}^{\delta}\V^{\tilde\beta_\nu}\big)
$$
is a local subsystem of $\X$, defined over $\Q$. By assumption
both must be equal. One obtains
\begin{equation}\label{isom}
\psi: \X \> \simeq >> \bigoplus_{\nu=1}^{\delta}\big(
\bigotimes_{(\nu-1)\delta'+1}^{\nu\delta'} \L_i \big)\otimes
{\T'}^{{\tilde\beta_\nu}}.
\end{equation}
Let us show next, that $\T'$ is a trivial local system. The
$\bar\Q$ isomorphism in (\ref{isom}) induces an isomorphism
$$
\END(\X) \> \simeq >>
\END\big(\bigoplus_{\nu=1}^{\delta}\V^{\tilde\beta_\nu}\big).
$$
Since $\beta \in {\rm Gal}(\bar\Q/\Q)$ permutes the direct factors
$\V^{\tilde\beta_\nu}$ of $\X$,
$$
\bigoplus_{\nu=1}^{\delta}\END(\V^{\tilde\beta_\nu})
$$
is a local subsystem, defined over $\Q$. So $\phi^{-1}$ induces
an embedding
$$
\phi': \bigoplus_{\nu=1}^{\delta}
\END(\L_{(\nu-1)\delta'+1})\otimes \cdots
\otimes\END(\L_{\nu\delta'})\otimes\END({\T'}^{\tilde\beta_\nu})
\>>> \END(\X),
$$
Writing $\END(\L_i)=\C\oplus\END_0(\L_i)$ we obtain a
decomposition of the left hand side in direct factors, all of the
form
$$
\END_0(\L_{j_1})\otimes \cdots \otimes
\END_0(\L_{j_\ell})\otimes\END({\T'}^{\tilde\beta_\nu}),
$$
for some $(\nu-1)\delta'+1 \leq j_1 < \cdots < j_\ell \leq \nu
\delta'$.

The only ones, without any $\END_0(\L_i)$ are the
$\END({\T'}^{\tilde\beta_\nu})$. We claim that
$$
\phi'(\bigoplus_{\nu=1}^{\delta}\END({\T'}^{\tilde\beta_\nu}))^{\beta}
=\phi'(\bigoplus_{\nu=1}^{\delta}\END({\T'}^{\tilde\beta_\nu})),
$$
for all $\beta\in {\rm Gal} (\bar\Q/\Q)$. Otherwise, we would get
a non-zero projection from
$\phi'(\bigoplus_{\nu=1}^{\delta}\END({\T'}^{\tilde\beta_\nu}))$
to an irreducible local system $\E$, containing at least one of
the $\END_0(\L_i)$. By construction, there exists an $\beta_i\in
{\rm Gal}(\bar\Q/\Q),$ such that $\L_i^{\beta_i}$ has a maximal
Higgs field. Hence $\E^{\beta_i}$ has a maximal Higgs field.

Applying $\beta_i$ we obtain a non-zero map
$$
\phi'(\bigoplus_{\nu=1}^{\delta}\END({\T'}^{\tilde\beta_\nu}))^{\beta\beta_i}
\>>> \E^{\beta_i}.$$ The right hand side has a maximal Higgs
field induced by the one on $\END_0(\L_i^{\beta_i})$, whereas the
left hand side is unitary, a contradiction.

So
$\phi'(\bigoplus_{\nu=1}^{\delta}\END({\T'}^{\tilde\beta_\nu}))$
is ${\rm Gal}(\bar\Q/\Q)$ invariant, hence a unitary local system
admitting a $\Z-$structure. This implies that
$\phi'(\bigoplus_{\nu=1}^{\delta}\END({\T'}^{\tilde\beta_\nu}))$
is trivial, after replacing $Y$ by a finite \'etale cover. So the
same holds true for $\END(\T')$, hence for $\T'$ as well. Let us
write $\T'=\C^{\oplus \epsilon'}$. Hence for some $\epsilon'$ one
has the decomposition (\ref{dec}), and a) holds true by
construction.

Recall that the local system $\L$ is defined over $F(\sqrt{a})$
for $a$ as in \ref{ramified}. Hence $\L_i$ is defined over
$\sigma_i(F)(\sqrt{\sigma_i(a)})$, and $\L_1\otimes \cdots
\otimes \L_{\delta'}$ is defined over the compositum $F'$ of those
fields, for $i=1, \ldots, \delta'$.

By \ref{quatsplit} $\Cor_{F/L}A$ can only split if $L=\Q$ and if
$[F:\Q]$ is odd. Let us write $L'=\Q$ in this case. Otherwise it
splits over the subfield $L'=L(\sqrt{b})$ of $F'$, where $b$ is
given in \ref{quatsplit}, b). In both cases one finds
$$
(\Cor_{F/L}A)\otimes_LL' \simeq M(2^{\delta'},L')
$$
and correspondingly $\L_1\otimes \cdots \otimes \L_{\delta'}$ is
defined over $L'$.

If $L'=\Q$, this is a local subsystem of $\X_\Q$. Since it is a
$\Q$ variation of Hodge structures, and since we assumed $\X_\Q$
to be irreducible, both coincide.

If $L'\neq L$ consider the $L'$ local subsystem $(\L_1\otimes
\cdots \otimes \L_{\delta'})_{L'}$ of $\V_{L'}$. For $\iota$ as
in d),
$$
\V' = (\L_1\otimes \cdots \otimes \L_{\delta'}) \oplus
(\L_1\otimes \cdots \otimes \L_{\delta'})^\iota,
$$
is a local subsystem of $\V$, defined over $L$, and of rank $
2^{\delta'+1}$. Then
$$
\bigoplus_{nu=1}^\delta\psi^{-1}({\V'})^{\tilde\beta_\nu}
$$
is a local subsystem of rank $\delta \cdot 2^{\delta'+1}$ of $\X$,
defined over $\Q$. It is also a sub variation of Hodge
structures. Since we assumed $\X_\Q$ to be irreducible, both must
coincide and $\epsilon'$ is equal to two.

It remains to verify e). Assume that $\M$ is a direct factor of
$\L_1\otimes \cdots \otimes \L_{\delta'}$. By \ref{fieldofdef2}
we may assume that $\M$ is defined over $\bar\Q$.

\ref{width}, i), implies that $\M$ has a maximal Higgs field. By
\ref{semistableii} $\M=\L'\otimes \T'_1$, and replacing $Y$ by an
\'etale covering we may assume that $\L'=\L_1=\L$, and that $\T'$
is a direct factor of $\L_2\otimes \cdots \otimes \L_{\delta'}$.
Using the notations introduced in \ref{notations}, let
$\tilde\sigma_i\in {\rm Gal}(\bar\Q/L)$ be an extension of
$\sigma_i$, for $i=1,\ldots,\delta'$. For those $i$ by
\ref{quatsplit2}
$$
(\L_1\otimes \cdots \otimes \L_{\delta'})^{\tilde\sigma_i}
$$
has again a maximal Higgs field. Applying \ref{width}, i),
one obtains the same for
$$
\M^{\tilde\sigma_i}=\L_1^{\tilde\sigma_i}\otimes{\T'_1}^{\tilde\sigma_i}.
$$
For $i=2$, the first factor is unitary, hence the second has
again a maximal Higgs field. \ref{addendum} tell us, that
replacing $Y$ again by some \'etale covering,
$$
{\T'_1}^{\tilde\sigma_i}=\L\otimes \T'',
$$
hence for $\T'_2={\T''}^{\tilde\sigma_i^{-1}}$
$$
\M=\L_1\otimes \L_2 \otimes \T'_2.
$$
Repeating this construction one finds
$$
\M=\L_1\otimes \cdots\otimes \L_{\delta'} \otimes \T'_{\delta'},
$$
necessarily with $\T'_{\delta'}=\C$.
\end{proof}

\begin{proposition}\label{rigidity}
Let $f:X\to Y$ be a family of abelian varieties with
general fibre $X_\eta$, and reaching the Arakelov bound. Then
\begin{enumerate}
\item[i.] For a generic fibre $X_\eta$ of $f$
$$
\End(X_\eta)\otimes \Q \simeq \End_Y(X)\otimes \Q \simeq
\End(R^1f_*\Q_{X})^{0,0}.
$$
\item[ii.] If $R^1f_*\C_{X}$ has no unitary part then
\begin{enumerate}
\item[a.] $\End(R^1f_*\Q_{X})^{0,0}=\End(R^1f_*\Q_{X})$.
\item[b.] If $X_\eta$ is geometrically simple, $R^1f_*\Q_{X}$ is irreducible.
\item[c.] $f:X\to Y$ is rigid, i.e. the morphism from $Y$
to the moduli scheme of polarized abelian varieties has no
non-trivial deformation.
\end{enumerate}
\end{enumerate}
\end{proposition}
\begin{proof}
i) is a special case of \cite{Del}, 4.4.6..

If $R^1f_*\C_X$ has no unitary part, for $\V=R^1f_*\C_X$
\ref{split3} gives a decomposition
$\END(\V)=\W\oplus \U$ where $\W$ has a maximal
Higgs field, and where $\U$ is concentrated in bidegree $0,0$. Since
\ref{width}, a), implies that $\W$ has no global section, one gets a).

For $X_\eta$ geometrically simple $\End(X_\eta)\otimes
\Q=\End(R^1f_*\Q_{X})^{0,0}$ is a skew field, hence a) implies that
$R^1f_*\Q_{X}$ is irreducible.

ii), c), follows from \cite{Fal} (see also
\cite{Sai}).
\end{proof}
\begin{proposition}\label{rigidity2}
Let $f:X\to Y$  be a family of abelian varieties, with a geometrically simple
generic fibre
$X_\eta$ and reaching the Arakelov bound. Assume that (replacing $Y$ by
an \'etale covering, if needed) one has the decomposition
(\ref{dec}) in \ref{quatsplit3}. Then $R^1f_*\C_{X}$ has no
unitary part if and only if
\begin{equation}\label{zerozero}
\End(R^1f_*\Q_{X})^{0,0}=\End(R^1f_*\Q_{X}).
\end{equation}
\end{proposition}
\begin{proof}
By \ref{rigidity} ii), a) and b), if $\X_\Q=R^1f_*\Q_X$ has no unitary part,
$\X_\Q$ is irreducible, and (\ref{zerozero}) holds true.

If on the other hand, $R^1f_*\C_X$ has a unitary part, the same
holds true for $\X$. Let us write again $\U_1$ for the unitary part
of $\X$. So the field $L$ in \ref{assumption2} can not be $\Q$.
Recall that the Higgs field of $\U_1$ splits in two components,
one of bidegree $1,0$, the other of bidegree $0,1$, both with a
trivial Higgs field. Correspondingly $\U_1$ is the direct sum of
two subsystems, say $\U^{1,0}$ and $\U^{0,1}$.

By \ref{quatsplit3} $\L_1\otimes \cdots \otimes \L_{\delta'}$
is an irreducible $\C$ local system. Let us choose one element of
$\C^{\epsilon'}$ and the corresponding local subsystem
$\M=\psi^{-1}(L_1\otimes \cdots \otimes \L_{\delta'})$ of $\X$.
There exists some $\beta\in
{\rm Gal}(\bar\Q/\Q)$ with $\M^\beta$ and $\bar\M^\beta$ unitary.
Replacing $\M$ by $\bar\M$, if necessary we may assume that
$\M^\beta$ lies in $\U^{1,0}$ and $\bar\M^\beta$ in $\U^{0,1}$.
Then
$$
\M^\beta \otimes \bar\M^{\beta^\vee} \subset
\U^{1,0}\otimes{\U^{0,1}}^\vee \subset \END(R^1f_*\C_X)^{1,-1}.
$$
In \ref{quatsplit}, v), we have seen that $\L_i\simeq \bar\L_i$
for all $i$. Hence
$$
\bar\M\simeq \bar\L_1 \otimes \cdots \otimes \bar\L_{\delta'}\simeq \M,
$$
and $\M^\beta$ and $\bar\M^\beta$ are isomorphic. One obtains
$\End(R^1f_*\C_X)^{1,-1}\neq 0$.
\end{proof}

\begin{proof}[Proof of \ref{shimura}]
Replacing $Y$ by an \'etale covering, we may assume that
$Rf_*\C_X$ has no unitary part as all. \ref{semistableii} provides
us with a local system $\L$, independent of all choices, again
after replacing $Y$ by some \'etale covering.

Hence it is sufficient to consider the case that the generic fibre
of $f:X\to Y$ is geometrically simple. By \ref{rigidity}, iv), the local system
$\X_\Q=R^1f_*\Q_X$ is irreducible. In \ref{quatsplit3} the non
existence of a unitary part implies that $\delta=1$, hence $L=\Q$,
and
$$
\X=\V=(\L_1\otimes \cdots \otimes \L_d)^{\oplus \epsilon'}.
$$
For $\epsilon'=1$, the $\Q$ local system $\X_\Q$ is given by the
representation
$$
\eta:\pi_1(Y,*) \>>> D^*=(\Cor_{F/\Q}A)^* = \Gl(2^{d},\Q).
$$
By \ref{quotient} $\pi_1(Y,*)\to \Gamma=\eta(\pi_1(Y,*)$ is an
isomorphism and $Y=\sH/\Gamma$. Hence $\X_\Q$ is isomorphic to
the local system $\X_{A\Q}$ constructed in \ref{constrq}. In
particular, $d=[F:\Q]$ is odd, and by \ref{rigidity}, i),
\ref{rigidity2}, and \ref{constrq1}
$$
\End(X_\eta)=\End(\X_\Q)=\Q \mbox{ \ \ and \ \ }
H^0(Y,\X_\Q\otimes \X_\Q)=\Q.
$$
The second equality implies that the polarization of $\X_\Q$ is
unique, up to multiplication with constants, hence $\X_\Q$ and
$\X_{A\Q}$ are isomorphic as polarized variations of Hodge
structures. For some $\Z$ structure on $\X_{A\Q}$ we constructed
in \ref{constrq2} a smooth family of abelian varieties $X_A \to
Y$, and this family is isogenous to $f:X\to Y$. Both satisfy the
properties, stated in Example \ref{mumfordshimura}, i).

For $\epsilon'=2$ and for $b$ as in \ref{quatsplit}, $\X_\Q$ is
given by
\begin{gather*}
\pi_1(Y,*) \>>> D^*=(\Cor_{F/\Q}A)^* \subset
(D\otimes_\Q\Q(\sqrt{b}))^*\\
= \Gl(2^d,\Q(\sqrt{b})) \subset \Gl(2^{d+1},\Q),
\end{gather*}
hence again $X_\Q$ is isomorphic to the local system $\X_{A\Q}$
constructed in \ref{constrq}.

By \ref{rigidity}, i), \ref{rigidity2} and \ref{constrq1}, i), one
finds that
$$
\End(X_\eta))=\End(\X_\Q)^{0,0}=\End(\X_\Q),
$$
is of dimension $4$.

For $b$ as in \ref{quatsplit}, consider the local system
$$
\L_{1\Q(\sqrt{b})}\otimes \cdots \otimes \L_{d\Q(\sqrt{b})}
$$
defined by the representation $\pi_1(Y,*) \to
\Gl(2^d,\Q(\sqrt{b}))$, together with a embedding into
$\X_{\Q(\sqrt{b})}$. Restricting the polarization, one obtains a
polarization $Q'$ on $\L_{1\Q(\sqrt{b})}\otimes \cdots \otimes
\L_{d\Q(\sqrt{b})}$, unique up to multiplication with constants.
Regarding this local system as a $\Q$ local system, the inclusion
$$
\Gl(2^d,\Q(\sqrt{b})) \subset \Gl(2^{d+1},\Q)
$$
defines an isomorphism
$$
\L_{1\Q(\sqrt{b})}\otimes \cdots \otimes \L_{d\Q(\sqrt{b})}\>>>
\X_\Q
$$
and the restriction of the polarization of $\X_\Q$ is the
composite of $Q'$ with the trace on $\Q(\sqrt{b})$. In
particular, the polarization is uniquely determined, and the
family $f:X\to Y$ is isogenous to the family $X_A\to Y_A=Y$
constructed in \ref{constrq2}.

Since, up to a shift in the bidegrees,
$$
\R^2f_*\Q_X=\bigwedge^2\X_\Q
$$
is a sub variation of Hodge structures of $\END(\X_\Q)$ one
obtains the first equality in
$$
\dim(H^0(Y,R^2f_*\Q)^{1,1})=\dim(H^0(Y,R^2f_*\Q))=\left\{
\begin{array}{ll}
3 & \mbox{for } d \mbox{ odd}\\
1 & \mbox{for } d \mbox{ even} \end{array}\right. ,
$$
whereas the second one has been verified in \ref{constrq1}, ii).
$\dim(H^0(Y,R^2f_*\Q)^{1,1})$ is the Picard number of a general
fibre of $f:X\to Y$. In fact, the Neron-Severi group of a general
fibre is invariant under the special Mumford-Tate group of the
fibre, hence by \ref{hg}, a), it coincides with
$\dim(H^0(Y,R^2f_*\Q)^{1,1})$.

Looking to the list of possible Picard numbers and to the
structure of the corresponding endomorphism algebras for simple
abelian varieties (for example in \cite{BL}, p. 141), one finds
that $\End(X_\eta)\otimes \Q$ is a quaternion algebra over $\Q$,
totally indefinite for $d$ odd, and totally definite otherwise.
Hence $f:X\to Y$ satisfies the properties stated in Example
\ref{mumfordshimura}, ii).
\end{proof}

\begin{proof}[Proof of \ref{shimura2}]
Again we may assume that $R^1f_*\C_X$ has no non trivial unitary
subbundle defined over $\Q$. Let $\V\oplus\U_1$ be the
decomposition of $R^1f_*\C_X$ in a part with a maximal Higgs
field and a unitary bundle. By \ref{semistableii} one can write
$\V=\L\otimes \T$, where after replacing $Y$ by a finite
covering, $\L$ only depends on $Y$. If $h:Z\to Y$ is a sub family
of $f:X\to Y$ with a geometrically simple generic fibre, then repeating this
construction with $g$ instead of $f$, we obtain the same local
system $\L$, hence by \ref{takeuchi} the same quaternion algebra
$A$. Hence we may assume that $f:X \to Y$ has a geometrically simple generic fibre,
and we have to show, that $f:X\to Y$ is one of the families in Example
\ref{general}.

By \cite{Del}, \S 4, $R^1f_*\Q_X$ is a direct sum of the same
irreducible $\Q$ local system $\X_\Q$. From \ref{semistableii} and
\ref{takeuchi} we obtain $\L$ and a quaternion algebra $A$,
defined over a totally real number field $F$. By \ref{quatsplit2},
$\X$ contains a local system $\V$, defined over a subfield $L$ of
$F$. which satisfies the conditions stated there. By
\ref{quatsplit3}, for $b$ as in \ref{corestriction}, $\V$ is given
by the representation $\pi_1(Y,*)\to \Gl(2^{\delta'+1},L)$
induced by
$$
\pi_1(Y,*) \>>> D_L=\Cor_{F/L}A \subset
D_L\otimes_LL(\sqrt{b})=M(2^{\delta'},L(\sqrt{b})) \subset
M(2^{\delta'+1},L),
$$
hence it is isomorphic to the local system in \ref{constrl}. Then
the decomposition of $\X$ in direct factors in \ref{quatsplit3}
coincides with the one in \ref{constrl}, and $f:X \to Y$ is one
of the families in Example \ref{general}.

In iii), the condition b) implies a) and vice versa. On the other
hand, $L_i=\Q$ if and only if $R^1h_{i*}\C_{Z_i}$ has no unitary
part, which by \ref{rigidity2} is equivalent to c).
\end{proof}

\section{Families of curves and Jacobians}\label{sectjacobians}

Let us shortly discuss the relation between Theorems
\ref{geomsplit} and \ref{shimura} and the number of singular
fibres for semi-stable families of curves.

Let $Y$ be a curve, let $h: \sC\to Y$ be a semi-stable
non-isotrivial family of curves of genus $g>1$, smooth over $V$,
and let $f:J(\sC/Y)\to Y$ be a compactification of the Neron model
of the Jacobian of $h^{-1}(V)\to V$. Let us write $S$ for the
points in $Y-V$ with $f^{-1}(y)$ singular, and $\Upsilon$ for the
other points in $Y\setminus V$, i.e. for the points $y$ with
$h^{-1}(y)$ singular but $f^{-1}(y)$ smooth. Let $g(Y)$ be the
genus of $Y$ and $U=Y\setminus S$.

The Arakelov inequality for non-isotrivial families of curves says
that
\begin{equation}\label{aracurves}
0 < 2\cdot \deg(F^{1,0}) \leq g_0 \cdot (2\cdot g(Y) - 2 + \# S +
\# \Upsilon),
\end{equation}
whereas the Arakelov inequality for $f:J(\sC/Y)\to Y$ gives the
stronger bound
\begin{equation}\label{arajac}
0 < 2\cdot \deg(F^{1,0}) \leq g_0 \cdot (2\cdot g(Y) - 2 + \# S ).
\end{equation}
Hence for a family of curves, the right hand side of
(\ref{aracurves}) can only be an equality, if $\Upsilon$ is empty.
On the other hand, if both, $S$ and $\Upsilon$ are empty, the
Miyaoka-Yau inequality for the smooth surface $\sC$ implies that
\begin{equation*}\label{my}
\deg(h_*\omega_{\sC/Y}) \leq \frac{g-1}{6} (2\cdot g(Y) - 2).
\end{equation*}
Hence if $h:\sC \to Y$ is smooth and if $h_*\C_{\sC}$ has no
unitary part, the inequalities (\ref{aracurves}) and
(\ref{arajac}) have both to be strict.\\

Let us consider the case $g(Y)=0$, i.e. families of curves over
$\P^1$. S.-L. Tan \cite{Tan} has shown that $h:\sC\to \P^1$ must
have at least $5$ singular fibres, hence that $\# S + \# \Upsilon
\geq 5,$ and (\ref{aracurves}) is strict in this special case.

Moreover, he and Beauville \cite{Bea1} gave examples of families
with exactly $5$ singular fibres for all $g>1$. In those examples
one has $\Upsilon = \emptyset$.

On the other hand, (\ref{arajac}) implies that $\# S \geq 4$. For
$\# S = 4$, the family $f:J(\sC/Y) \to Y$ reaches the Arakelov
bound, hence by \ref{geomsplit} it is isogenous to a product of a
constant abelian variety with a product of modular elliptic
curves, again with $4$ singular fibres. By \cite{Bea2} there are
just $6$ types of such families, among them the universal family
$E(3)\to X(3)$ of elliptic curves with a level $3$-structures.

Being optimistic one could hope, that those families can not occur
as families of Jacobians, hence that there is no family of curves
$h:\sC\to P^1$ with $\# S = 4$. However, a counterexample has been
constructed in \cite{Kan1}.

\begin{example}\label{elldiff}
Let $B$ be a fixed elliptic curve, defined over $\C$. Consider the
Hurwitz functor $\sH_{B,N}$ defined in \cite{Kan1}, i.e. the
functor from the category of complex schemes to the category of
sets with
\begin{multline*}
\sH_{B,N}(T)=\{f: C\>>> B\times T; \ f \mbox{ is a normalized
covering of degree }N \\
\mbox{ and }C\mbox{ a smooth family of curves of genus } 2 \mbox{
over } T \}.
\end{multline*}
The main result of \cite{Kan1} says that for $N \geq 3$ this
functor is represented by an open subscheme $V=H_{B,N}$ of the
modular curve $X(N)$ parameterizing elliptic curves with a level
$N$-structure.

The universal curve $\sC \to H_{B,N}$ extends to a
semi-stable curve $\sC\to X(N)$ whose Jacobian is isogenous to
$B\times E(N)$. Hence writing $S$ for the cusps, $J(\sC/X(N))$ is
smooth outside of $S$, whereas $\sC\to X(N)$ has singular
semi-stable fibres outside of $H_{B,N}$. Theorem 6.2 in
\cite{Kan1} gives an explicit formula for the number of points in
$\Upsilon=X(N)\setminus (H_{B,N}\cup S)$.

Evaluating this formula for $N=3$ one finds $\# \Upsilon = 3$. For
$N=3$ the modular curve $X(3)$ is isomorphic to $\P^1$ with $4$
cusps. So the number of singular fibres is $4$ for $J(\sC/\P^1)
\to \P^1$ and $7$ for $\sC \to \P^1$.
\end{example}

We do not know whether similar examples exist for $g>2$. For
$g>7$ the constant part $B$ in Theorem \ref{geomsplit} can not be
of codimension one. In fact, the irregularity $q(\sC)$ of the
total space of a family of curves of genus $g$ over a curve of
genus $q$ satisfies by \cite{Xia}, p. 461, the inequality
$$
q(\sC) \leq \frac{5\cdot g +1}{6} + g(Y).
$$
If $J(\sC/Y) \to Y$ reaches the Arakelov bound, hence if it is
isogenous to a product
$$ B\times E\times_{Y}\cdots\times_{Y} E ,$$
one finds $$\dim(B)\leq\frac{5\cdot g +1}{6}.$$

As explained in \cite{ES} it is not known, whether for $g\gg 2$
there are any curves $C$ over $\C$ whose Jacobian is isogenous to
the product of elliptic curves. Here we are even asking for families
of curves whose Jacobian is isogenous to the product of the same
non-isotrivial family of elliptic curve, up to a constant factor.

For the smooth families of abelian varieties, considered in
\ref{shimura} or \ref{shimura2} we do not know of any example,
where such a family is a family of Jacobians.

\bibliographystyle{plain}

\end{document}